\documentclass[12pt,leqno]{article}

\usepackage{amsfonts}
\usepackage{amsmath}
\usepackage{amssymb}

\newcommand{\rarrow}[1]{{\buildrel #1 \over \longrightarrow}}

\def\Im{{\rm Im}}

\def\Z{{\mathbb Z}}
\def\C{{\mathbb C}}
\def\R{{\mathbb R}}
\def\H{{\mathbb H}}

\def\S{{\mathbb S}}
\def\P{{\rm P}}

\date      \kill

\renewcommand{\k}{\hspace{-2mm}{\bf .}\hspace{5mm}}

\title{\sc Gottlieb groups of spheres}
\author{Marek Golasi\'nski and Juno Mukai}

\newtheorem{thm}{Theorem}[section]
\newtheorem{cor}[thm]{Corollary}
\newtheorem{prop}[thm]{Proposition}
\newtheorem{lem}[thm]{Lemma}
\newtheorem{rem}[thm]{Remark}

\newtheorem{exam}[thm]{Example}
\newtheorem{con}[thm]{Conjecture}

\begin{document}

\maketitle

{\bf Abstract.} This paper takes up the systematic study of the Gottlieb groups
$G_{n+k}(\S^n)$ of spheres for $k\le 13$ by means of the classical homotopy
theory methods. The groups $G_{n+k}(\S^n)$ for $k\le 7$ and $k=10,12,13$ are fully
determined. Partial results on $G_{n+k}(\S^n)$ for $k=8,9,11$ are presented as well.
We also show that $[\iota_n,\eta^2_n\sigma_{n+2}]=0$ if $n=2^i-7$ for $i\ge 4$.

\footnote{
The first author was partially supported by Akio Sat$\bar{\rm o}$, Shinshu University 101501.
The second author was partially supported by Grant-in-Aid for Scientific Research (No.\ 15540067 (c), (2)), Japan
Society for the Promotion of Science and the Faculty of Mathematics and Computer Science in Toru\'n.}

\section*{Introduction}
The Gottlieb groups $G_k(X)$ of a pointed space
$X$ have been defined by Gottlieb in \cite{G} and \cite{G1}; first $G_1(X)$ and then
$G_k(X)$ for all $k\ge 1$.  The higher Gottlieb
groups $G_k(X)$ are related in \cite{G1} and \cite{G2} to the existence of
sectioning fibrations with fiber $X$. For instance, if $G_k(X)$ is trivial
then there is a homotopy section for every fibration over the $(k+1)$-sphere
$\S^{k+1}$, with fiber $X$.
\par This paper grew out of our attempt to develop techniques in
calculating $G_{n+k}(\S^n)$ for $k\le 13$ and any $n\ge 1$.
The composition methods developed by Toda \cite{T} are the main
tools used in the paper. Our calculations also deeply depend on the results
of \cite{H-M}, \cite{K} and \cite{M2}.
\par Section 1 serves as backgrounds to the rest of the paper.
Write $\iota_n$ for the homotopy class of the identity map of
$\S^n$. Then, the homomorphism
$$P' : \pi_k(\S^n)\longrightarrow\pi_{k+n-1}(\S^n)$$
defined by $P'(\alpha)=[\iota_n,\alpha]$
for $\alpha\in\pi_k(\S^n)$ \cite{Hilton1} leads to
the formula $G_k(\S^n)=\mbox{ker}\,P'$, where $[-,-]$
terms the standard Whitehead product. So, our main task is to consult first
\cite{Hilton1}, \cite{Hil}, \cite{M1}, \cite{M2}, \cite{Thomeier} and \cite{T}
about the order of $[\iota_n,\alpha]$ and
then to determine some Whitehead products in unsettled cases as well.
In the light of Serre's result \cite{Serre}, the $p$-primary component of
$G_{2m+k}(\S^{2m})$ vanishes for any odd prime $p$, if $2m\ge k+1$ (Proposition \ref{Gnp0}).
\par Let $EX$ be the suspension of a space $X$ and denote by
$E: \pi_k(X)\to\pi_{k+1}(EX)$ the suspension map. Write
$\eta_2\in\pi_3(\S^2)$, $\nu_4\in\pi_7(\S^4)$ and $\sigma_8
\in\pi_{15}(\S^8)$ for the Hopf maps, respectively. We set
$\eta_n=E^{n-2}\eta_2\in\pi_{n+1}(\S^n)$ for $n\ge 2$,
$\nu_n=E^{n-4}\nu_4\in\pi_{n+3}(\S^n)$ for $n\ge 4$ and
$\sigma_n=E^{n-8}\sigma_8\in\pi_{n+7}(\S^n)$ for $n\ge 8$.
Write $\eta^2_n=\eta_n\circ\eta_{n+1}$, $\nu^2_n=\nu_n\circ\nu_{n+3}$ and
$\sigma^2_n=\sigma_n\circ\sigma_{n+7}$.
Section 2 is a description of $G_{n+k}(\S^n)$ for $k\le 7$. To reach that
for $G_{n+6}(\S^n)$, we make use of  Theorem \ref{neq} partially extending the
result of \cite{K-M}: {\em $[\iota_n, \nu^2_n] = 0$ if and only if
$n\equiv 4,5,7\ (\bmod\ 8)$, $n = 2^i - 3$ or $n = 2^i - 5$ for $i\ge 4$;}
for the proof of which Section 3 and Section 4 are devoted.
\par Section 5 devotes to proving Mahowald's claim:
$[\iota_n,\sigma_n]\ne 0$ for $n\equiv 7\ (\bmod\ 16)\ge 23$.
\par Section 6 takes up computations of $G_{n+k}(\S^n)$ for
$k=10,12,13$ and partial ones of $G_{n+k}(\mathbb{S}^n)$ for $k=8,9,11$.
In a repeated use of \cite{M2}, we have found out the triviality of
the Whitehead product \mbox{\cite{M5}}:
$$
[\iota_n,\eta^2_n\sigma_{n+2}]=0, \ \mbox{if}\;
n=2^i-7\ (i\ge 4),$$
which corrects thereby \cite{M2} for $n=2^i-7$.
\par The authors thank Professor M.\ Mimura for suggesting the problem and
fruitful conversations, Professors H.\ Ishimoto, I.\ Madsen, M.\ Mahowald,
Y.\ Nomura and N.\ Oda for helpful informations on the orders of the Whitehead products
$[\iota_n,\nu^2_n]$ and $[\iota_n,\sigma_n]$.
\par The authors are also very grateful to Professor H.\ Toda for informing
the order of the Whitehead product $[\iota_{2n},[\iota_{2n},\alpha_1(2n)]]$ for $n\ge
2$ (Proposition \ref{Toda}), where $\alpha_1(2n)$ is a generator of the
$3$-primary component of $\pi_{2n+3}(\mathbb{S}^{2n})$.

\section{Preliminaries on Gottlieb groups}\setcounter{equation}{0}
\par Throughout this paper, spaces, maps and homotopies are based. We use the standard terminology
and notations from the homotopy theory, mainly from \cite{T}. We do not distinguish
between a map and its homotopy class.
\par Let $X$ be a connected space.
The $k$-th {\em Gottlieb group} $G_k(X)$ of $X$ is the subgroup of
the $k$-th homotopy group $\pi_k(X)$ consisting of all elements which can be
represented by a map $f : \S^k\to X$ such that $f\vee\mbox{id}_X : \S^k\vee
X\to X$ extends (up to homotopy) to a map $F : \S^k\times X\to X$.
Define $P_k(X)$ to be the set of elements of $\pi_k(X)$ whose Whitehead product
with all elements of all homotopy groups is zero. It turns out that $P_k(X)$
forms a subgroup of $\pi_k(X)$ and, by \cite[Proposition 2.3]{G1}, $G_k(X)\subseteq P_k(X)$.
Recall that $X$ is said to be a $G$-{\em space} ({\em resp.} $W$-{\em space}) if
$\pi_k(X)=G_k(X)$ ({\em resp.} $\pi_k(X)=P_k(X)$) for all $k$.
\par Given $\alpha\in\pi_k(\S^n)$ for $k\ge 1$, we deduce that
$\alpha\in G_k(\S^n)$  if and only if $[\iota_n,\alpha]=0$.
In other words, consider the map $$P' : \pi_k(\S^n)\longrightarrow
\pi_{k+n-1}(\S^n)$$ defined by $P'(\alpha)=[\iota_n,\alpha]$ for $\alpha\in\pi_k(\S^n)$.
Then, this leads to the formula $$G_k(\S^n)=\mbox{ker}\,P'.$$
\par Write now $\sharp$ for the order of a group or its any element.
Then, from the above interpretation of Gottlieb groups of spheres, we obtain
\begin{lem}\label{lead}
If $\pi_k(\S^n)$ is a cyclic group for some $k\ge 1$ with a gene\-ra\-tor
$\alpha$ then $G_k(\S^n)=(\sharp[\iota_n,\alpha])\pi_k(\S^n)$.
\end{lem}
Since $\S^n$ is an H-space for $n=3,7$, we have $$G_k(\S^n) = \pi_k(\S^n)\ \mbox{for}\ k\ge 1,
\ \mbox{if} \ n =3,7.$$
\par We recall the following result from \cite{Hil} and \cite{WG} needed in the sequel.
\begin{lem} \label{pro}
\mbox{\em (1)} If $\xi\in\pi_m(X)$, $\eta\in\pi_n(X)$, $\alpha\in\pi_k(\S^m)$,
$\beta\in\pi_l(\S^n)$ and if $[\xi,\eta]=0$ then
$[\xi\circ\alpha,\eta\circ\beta]=0.$

\mbox{\em (2)} Let $\alpha\in\pi_{k+1}(X)$, $\beta\in\pi_{l+1}(X)$,
$\gamma\in\pi_m(\S^k)$ and
$\delta\in\pi_n(\S^l)$.

Then $[\alpha\circ E\gamma,\beta\circ E\delta]=[\alpha,\beta]\circ E(\gamma\wedge\delta).$

\mbox{\em (3)} If $\alpha\in\pi_k(\S^2)$ and $\beta\in\pi_l(\S^2)$ then
$[\alpha,\beta]=0$ unless $k=l=2$.

\mbox{\em (4)} $[\beta,\alpha] = (-1)^{(k+1)(l+1)}[\alpha,\beta]$ for
$\alpha\in\pi_{k+1}(X)$ and $\beta\in\pi_{l+1}(X)$.

In particular, $2[\alpha,\alpha]=0$ for $\alpha\in\pi_n(X)$ if $n$ is odd.

\mbox{\em (5)} If $\alpha_1,\alpha_2\in\pi_{p+1}(X)$,
$\beta\in\pi_{q+1}(X)$ and $p\ge 1$, then $[\alpha_1+\alpha_2,\beta]=$

$[\alpha_1,\beta]+[\alpha_2,\beta]$ and $[\beta,\alpha_1+\alpha_2]=[\beta,\alpha_1]
+[\beta,\alpha_2]$.

\mbox{\em (6)} $E[\alpha,\beta]=0$ for $\alpha\in\pi_k(X)$ and $\beta\in\pi_l(X)$.

\mbox{\em (7)} $3[\alpha,[\alpha,\alpha]]=0$ for $\alpha\in\pi_{n+1}(X)$.
\end{lem}

\par Let $G_k(X;p)$ and $\pi_k(X;p)$ be the $p$-primary components of $G_k(X)$
and $\pi_k(X)$ for a prime $p$, respectively. But for $X=\S^n$, recall the
notation
from \cite{T}:
$$
\pi^n_k= \left\{ \begin{array}{ll}
\pi_n(\S^n), & \mbox{if $k=n$};\\
E^{-1}\pi_{2n}(\S^{n+1};2),& \mbox{if $k=2n-1$};\\
\pi_k(\S^n;2),& \mbox{if $k\not=n,2n-1.$}
\end{array}
\right.
$$
\par As it is well-known, $[\iota_n,\iota_n]=0$ if and only if $n=1,3,7$ and
$\sharp[\iota_n,\iota_n]=2$ for $n$ odd and $n\ne 1,3,7$, and it is infinite provided
$n$ is even. Thus, we have reproved the result \cite{G1} that $G_n(\S^n)=\pi_n(\S^n)\cong\Z$ for
$n=1,3,7$, $G_n(\S^n)=2\pi_n(\S^n)\cong 2\Z$ for $n$ odd and $n\not=1,3,7$, and
$G_n(\S^n)=0$ for $n$ even, where $\Z$ denotes the additive group of integers.
It is easily obtained that $G_k(\S^n)=P_k(\S^n)$ for all $k,n$ \cite[Theorem I.9]{La}.
In other words, on the level of spheres the class of $G$-spaces coincides with that
of $W$-spaces.
\par Let now $n$ be odd. Then, by Lemma \ref{pro}:(4), (5) and (7),
$[\iota_n, [\iota_n,\iota_n]]=0$.
Furthermore, by Lemma \ref{pro}:(1), (4)
and (5), $[2\iota_n,\iota_n]=0$ implies that $[\iota_n,2\alpha]=0$ for any
$\alpha\in\pi_k(\S^n)$, that is, $2\pi_k(\S^n)\subset G_k(\S^n)$ and thus $G_k(\S^n;p)=\pi_k(\S^n;p)$
for any odd prime $p$. In the light of \cite{Liulevicius}, we also know
%\addtocounter{equation}{1}
\begin{equation} \label{Lv}
\sharp[\iota_{2n}, [\iota_{2n}, \iota_{2n}]]=3, \ \mbox{if} \ n\geq 2.
\end{equation}
Whence, Lemma \ref{pro} and (\ref{Lv}) yield the results proved partially in \cite{GG}.
\begin{cor}\k \label{CC}
{\em (1)} $3[\iota_n,\iota_n]\in G_{2n-1}(\S^n)$. In particular, $a[\iota_n,\iota_n]\in$

$G_{2n-1}(\S^n)$ for $n$ even if and only if $a\equiv 0\ (\bmod\ 3)$.

{\em (2)} If $k\ge 3$ then $G_k(\S^2)=\pi_k(\S^2)$.

{\em (3)} If $n$ is odd and $n\not=1,3,7$ then $2\pi_k(\S^n)\subset G_k(\S^n)$.
In particular,

$G_k(\S^n;p)=\pi_k(\S^n;p)$ for any odd prime $p$ and $k\ge 1$.

{\em (4)} $G_k(\S^n)=\pi_k(\S^n)$ provided that
$E : \pi_{k+n-1}(\S^n)\to\pi_{k+n}(\S^{n+1})$ is

a monomorphism.
\end{cor}

\par We note that $P'$ and the homomorphism
$$P :\pi_{k+n+1}(\mathbb{S}^{2n+1})
\longrightarrow\pi_{k+n-1}(\mathbb{S}^n)\ (k\le 2n-2)$$
in the EHP sequence are related as follows:
$$P'=P\circ E^{n+1}\;\mbox{for}\;k\le 2n-2.$$
\par Let $SO(n)$ be the rotation group and
$J:\pi_k(SO(n))\to\pi_{n+k}(\S^n)$ be the $J$-homomorphism and
$\Delta: \pi_k(\S^n)\to\pi_{k-1}(SO(n))$ the connecting map
associated with the fibration $SO(n+1)\stackrel{SO(n)}{\longrightarrow}\S^n$. We recall
\begin{equation} \label{P'}
P'=J\circ\Delta
\end{equation}
and so, %by the definition and (\ref{P'}), we obtain
\begin{equation} \label{JDel} \mbox{ker} \{\Delta: \pi_k(\S^n)\to\pi_{k-1}(SO(n))\}\subset G_k(\S^n).
\end{equation}
\par In virtue of \cite[Chapter IV]{Serre} (\cite[(13.1)]{T}), Serre's isomorphism
%\addtocounter{equation}{1}
\begin{equation} \label{Serre}
\pi_{i-1}(\S^{2m-1};p)\oplus\pi_i(\S^{4m-1};p)\cong\pi_i(\S^{2m};p)
\end{equation}
is given by the correspondence $(\alpha, \beta)\mapsto E\alpha
+ [\iota_{2m}, \iota_{2m}]\circ\beta$.
\par We show
\begin{prop}\k \label{Gnp0} Let $p$ be an odd prime and $n$ even.

\mbox{\em (1)} If $n\ge k+1$ then $G_{n+k}(\S^n;p)=0$.

\mbox{\em (2)} Suppose that $\pi_{n+k}(\S^n;p)$ is cyclic with a generator $\theta$,
$\pi_{n+k}(\S^n;p)=$

$E\pi_{n+k-1}(\S^{n-1};p)$ and that $E^{n-1}: \pi_{n+k}(\S^n;p)\to\pi_{2n+k-1}(\S^{2n-1};p)$

is an epimorphism. Then, $G_{n+k}(\S^n;p)=\{\sharp(\mbox{\em ker}\,E^{n-1})\theta\}$.
\end{prop}
{\bf Proof.}
By the Freudenthal suspension theorem, $E: \pi_{n+k}(\S^n)\to\pi_{n+k+1}(\S^{n+1})$
is an epimorphism, if $n=k+1$ and an isomorphism, if $n\ge k+2$. For the case
$n=k+1$, by the EHP sequence $$\pi_{2n-2}(\S^{n-1})\rarrow{E}\pi_{2n-1}(\S^n)\rarrow{H}\pi_{2n-1}(\S^{2n-1})$$
we get $$\pi_{2n-1}(\S^n)=E\pi_{2n-2}(\S^{n-1})\oplus\Z\{\alpha\},$$
where $\alpha=\eta_2,\nu_4,\sigma_8$ according as $n=2,4,8$ and $\alpha=[\iota_n,\iota_n]$ for $n$ otherwise.
So, by the EHP sequence
$$\pi_{2n+1}(\S^{2n+1})\rarrow{P}\pi_{2n-1}(\S^n)\rarrow{E}\pi_{2n}(\S^{n+1})$$
and by the Freudenthal suspension theorem,
$$\sharp\beta=\sharp(E^m\beta) \ \mbox{for any} \ \beta\in\pi_{n+k}(\S^n;p) \
\mbox{and} \ m\ge 1.$$
Hence, by Lemma \ref{pro}.(2) and (\ref{Serre}), we obtain
$$\sharp\beta=\sharp[\iota_n,\beta].$$
This leads to (1).
\par (2) is easily obtained from (\ref{Serre}) and the proof is complete.
\hfill$\square$

\par The notation $\pi_{n+m}(\S^n)=\{\alpha_n\}\ (\{\alpha(n)\}, \ resp.)$ means that there exist
some $k\ge 1$ and an element $\alpha_k\ (\alpha(k),\ resp.)\in\pi_{k+m}(\S^k)$ satisfying $\alpha_n=
E^{n-k}\alpha_k \ (\alpha(n)=E^{n-k}\alpha(k), \ resp.)$ for $n\ge k$. For the $p$-primary
component with any prime $p$, the notation is available.

%(Attention: We already used those notations and conventions. For example,
%$\pi_{n+1}(\S^n)=\{\eta_n\}$ for $n\ge 2$.)

Hereafter, we omit the reference \cite{T} unless otherwise stated.
Now, we know that $\pi_{n+3}(\S^n;3)=\{\alpha_1(n)\}\cong\Z_3$ and
$\pi_{n+7}(\S^n;3)=\{\alpha_2(n)\}\cong\Z_3$ for $n\ge 3$.
Write $\{-, -, -\}_n$ for the Toda bracket, where
$n\ge 0$ and $\{-, -, -\}=\{-, -, -\}_0$.
We recall that there exists the element $\beta_1(5)\in\pi_{15}(\S^5)$ satisfying
$\beta_1(5)\in\{\alpha_1(5),\alpha_1(8),\alpha_1(11)\}_1$,
$3\beta_1(5)=-\alpha_1(5)\alpha_2(8)$ and that $\pi_{n+10}(\S^n;3)=\{\beta_1(n)\}\cong\Z_9$ for $n=5,6$ and $\cong\Z_3$ for $n\ge 7$.

Let $\Omega^2\S^{2m+1}=\Omega(\Omega \S^{2m+1})$
be the double loop space of $\S^{2m+1}$ and $Q^{2m-1}_2
=\Omega(\Omega^2\S^{2m+1},\S^{2m-1})$ the homotopy fiber of the canonical
inclusion (the double suspension map) $i: \S^{2m-1}\to\Omega^2\S^{2m+1}$.
Then, the following result and its proof have been shown by Toda \cite{T4}.
%\addtocounter{thm}{1}
\begin{prop}\k \label{Toda}
Let $n$ be even and $n\ge 4$. Then $[\iota_n,[\iota_n,\alpha_1(n)]]\ne 0$
if and only if $n\ne 4$ and $n\equiv 1\ (\bmod\ 3)$.
\end{prop}
{\bf Proof.}
It is well-known that $[\iota_n,[\iota_n,\iota_n]]\in
E\pi_{3n-3}(\S^{n-1})$. So,
by the $(\bmod\ 3)$ EHP sequence \cite[(2.1.3), Proposition 2.1]{T3},
we have
$$
[\iota_n,[\iota_n,\iota_n]]=\pm EP(i(n-1)), \ \mbox{where} \ i(n-1) \ \mbox{is
a generator of} \ \pi_{3n-3}(Q^{n-1}_2;3).
$$
By the naturality \cite[(2.1.5)]{T3}, we obtain
$[\iota_n,[\iota_n,\alpha_1(n)]]
=\pm EP(a_1(n-1))$, where
$a_1(n-1)=i(n-1)\alpha_1(3n-3)\in\pi_{3n}(Q^{n-1}_2;3)$.
By \cite[Proposition 4.4]{T3}, $(\frac{n}{2}+1)a_1(n-1)=HP(i(n+1))$. So,
$P(a_1(n-1))=\pm PHP(i(n+1))=0$ if $n\not\equiv 1\ (\bmod\ 3)$. For the
case $n=4$, the assertion is trivial.

Next, assume that $n\ne 4$ and $n\equiv 1\ (\bmod\ 3)$. Then, by \cite[Theorem
10.3]{T2}, there exists an element $v\in\pi_{3n-4}(\S^{n-3})$ satisfying
$H(v)=b(n-5)$ and $E^2v=P(a_1(n-1))$, where $b(n-5)=i(n-5)\beta_1(3n-15)$.
Furthermore, by \cite[Proposition 5.3.(ii)]{T2}, we obtain $P(a_2(n-3))=3v$,
where $a_2(n-3)=i(n-3)\alpha_2(3n-9)$. So, by the $(\bmod\ 3)$ EHP sequence, we
have $P(a_1(n-1))\ne 0$. This implies the sufficient condition and
completes the
proof.
\hfill$\square$

\section{Gottlieb groups of spheres with stems for $k\le 7$}%{Gottlieb groups $G_{n+k}(\S^n)$ for $k\le 7$}
\setcounter{equation}{0}According to \cite{Hilton1}, \cite{Hil}, \cite{K-M}, \cite{M1}, \cite{Thomeier} and \cite{T},
we know the following results:
\begin{equation}[\iota_n, \eta_n]=0 \ \mbox{if and only if} \ n\equiv 3\, (\bmod\
4)\;\mbox{or}\;n=2,6;\label{weta}
\end{equation}
\begin{equation}[\iota_n, \eta^2_n]=0 \ \mbox{if and only if} \ n\equiv 2,3\, (\bmod\ 4)\;
\mbox{or} \; n=5.\label{weta2}
\end{equation}
Hence, Lemma \ref{lead} completely determines
$G_{n+k}(\S^n)$ for $k=1,2$.
\bigskip

\par We recall that $\pi^3_6=\{\nu'\}\cong\Z_4$. Write $\omega$ for a generator
of the $J$-image $J\pi_3(SO(3))=\pi_6(\S^3)\cong\Z_{12}$ satisfying $\omega=\nu'+\alpha_1(3)$.
We recall the relation $\pm[\iota_4,\iota_4]=2\nu_4-\Sigma\nu'$. By abuse of notation,
$\nu_n$ represents a generator of $\pi^n_{n+3}$ and $\pi_{n+3}(\S^n)$ for $n\ge 4$, respectively.
Then, $\pi_7(\S^4)=\{\nu_4,E\omega\}\cong\Z\oplus \Z_{12}$,
$\pi_{n+3}(\S^n)=\{\nu_n\}\cong\Z_{24}$ for $n\ge 5$.
By \cite{B-B}, $[\iota_4,\nu_4]=2\nu^2_4$.
In the light of Lemma \ref{pro}.(2) and the relation $\nu'\nu_6 = 0$, we obtain
$$[\iota_4, E\nu'] = [\iota_4, \iota_4]\circ E(\iota_3\wedge\nu')
= (2\nu_4 - E\nu')\circ 2\nu_7 = 4\nu^2_4.$$
So, we have $2E\nu' = \eta^3_4\in G_7(\S^4)$. Consequently, by Corollary \ref{CC}.(1) and
Proposition \ref{Gnp0}, $$G_7(\S^4) = \{3[\iota_4,\iota_4], 2E\nu'\}\cong 3\Z\oplus\Z_2.$$
In the light of \cite{K-M}, \cite{M1}, \cite{Thomeier}, \cite{T}, Corollary \ref{CC}.(3)
and Proposition \ref{Gnp0}, we know
the following:%\addtocounter{equation}{-3}
\begin{multline}\label{wnu}
\sharp[\iota_n,\nu_n] =\left\{\begin{array}{ll}
1,&\quad\mbox{if}\;n\equiv 7\; (\bmod\ 8), \; n = 2^i - 3\;
\mbox{for} \; i\ge 3;\\
2,&\quad\mbox{if} \; n\equiv 1,3,5\; (\bmod\ 8)\ge 9 \;
\mbox{and}\;n\neq 2^i - 3;\\
12,&\quad\mbox{if} \;n\equiv 2\;(\bmod\ 4)\ge 6, n = 12;\\
24,&\quad \mbox{if} \;n\equiv 0\;(\bmod\ 4)\ge 8 \ \mbox{unless} \; 12.
\end{array}
\right.
\end{multline}
Thus, Lemma \ref{lead} leads to a complete
description of $G_{n+3}(\S^n)$.

\bigskip

\par Now, we recall the following relations:
$$
\eta_n\nu_{n+1}=0 \ \mbox{for} \ n\ge 5 \ \mbox{and} \
\nu_n\eta_{n+3}=0 \ \mbox{for} \ n\ge 6.
$$
By the relation $\nu'\eta_6=\eta_3\nu_4$, we have
$[\iota_4,\nu_4\eta_7]=[\iota_4,(E\nu')\eta_7]
=[\iota_5,\nu_5\eta_8]=0$. Hence, by the group structures of
$\pi_{n+k}(\S^n)$ for $k = 4,5$, we get

\begin{prop}\k \label{prop3}
$G_{n+4}(\S^n)=\pi_{n+4}(\S^n)$; $G_{n+5}(\S^n)=\pi_{n+5}(\S^n)$ unless $n=6$
and $G_{11}(\S^6)=3\pi_{11}(\S^6)\cong 3\Z$.
\end{prop}

In the next two sections, we will prove the following result
partially extending that of \cite[Theorem 1.3]{K-M}.

\begin{thm}\k \label{neq}$[\iota_n, \nu^2_n] = 0$ if and only if
$n\equiv 4,5,7\ (\bmod\ 8)$, $n = 2^i - 3$ or $n = 2^i - 5$ for $i\ge 4$.
\end{thm}
We recall that $\pi_{10}(\S^4)=\{\nu^2_4,\alpha_1(4)\alpha_1(7),
\nu_4\alpha_1(7)\}\cong\Z_8\oplus(\Z_3)^2$.
By \cite{B-B} and the relation $\alpha_1(7)\alpha_1(10)=0$, we get that
$[\iota_4,\nu_4\alpha_1(7)]=[\iota_4,\alpha_1(4)\alpha_1(7)]=[\iota_4,\iota_4](\alpha_1(7)\alpha_1(10))=0$.
Recall also that $\pi^5_{12}=\{\sigma'''\}\cong\Z_2$, $\pi^6_{13}=\{\sigma''\}\cong\Z_4$,
$\pi^7_{14}=\{\sigma'\}\cong\Z_8$, where $E\sigma'''=2\sigma''$, $E\sigma''=2\sigma'$ and
$E^2\sigma'=2\sigma_9$. By \cite{B-B} and \cite{T}, we obtain
$$
[\iota_5,\sigma''']=[\iota_5,\iota_5]\circ E^4\sigma'''=0,
[\iota_6,\sigma'']=[\iota_6,\iota_6]\circ E^5\sigma''
= 4([\iota_6,\iota_6]\circ\sigma_{11})$$
and $2[\iota_6,\sigma'']\ne 0.$
We recall the relation $\pm[\iota_8,\iota_8]=2\sigma_8-E\sigma'$.
In $\pi^8_{22}=\Z_{16}\{\sigma^2_8\}\oplus\Z_8\{(E\sigma')\sigma_{15}\}\oplus\Z_4\{\kappa_8\}$,
we have $[\iota_8,E\sigma']=2[\iota_8,\iota_8]\sigma_{15}=\pm 2(2\sigma^2_8-(E\sigma')\sigma_{15})$
and in view of \cite{B-B}, we obtain $[\iota_8,\sigma_8]=[\iota_8,\iota_8]\circ\sigma_{15}=
\pm(2\sigma^2_8-(E\sigma')\sigma_{15})$.
We know that $\pi_{n+7}(\S^n;5)=\{\alpha'_1(n)\}$ $\cong\Z_5$ for $n\ge 3$.
Thus, by Corollary \ref{CC}, Proposition \ref{Gnp0} and
Theorem \ref{neq}, we obtain
\begin{prop}\k \label{ToMa}
$G_{n+6}(\S^n) = \pi_{n+6}(\S^n)$ if $n\equiv 4,5,7\ (\bmod\ 8)$,
$n = 2^i - 3$ or $n = 2^i - 5$ and $G_{n+6}(\S^n) = 0$ otherwise.
\par Furthermore, $G_{n+7}(\S^n)=0$ if $n=4,6$,
$G_{12}(\S^5)=\pi_{12}(\S^5)$, $G_{15}(\S^8)=\{3[\iota_8,\iota_8],4E\sigma'\}
\cong 3\Z\oplus\Z_2$.
\end{prop}

Let $H: \pi_k(\S^n)\to \pi_k(\S^{2n-1})$ be the Hopf homomorphism.
Then, by \cite{Ad} and \cite[Proposition 4.5]{Oda}, there exists an element
$\gamma\in\pi^{n-7}_{2n-8}$ satisfying%\addtocounter{equation}{-1}
\begin{equation} \label{des7}
[\iota_n,\iota_n]=E^7\gamma \ \mbox{and} \ H\gamma=\sigma_{2n-15}, \
\mbox{if} \ n\equiv 7\ (\bmod\ 16)\ge 23.
\end{equation}
According to Mahowald \cite{M5} and (\ref{des7}), we obtain
\begin{thm}\k \label{Mah}
$[\iota_n,\sigma_n]\ne 0$, if $n\equiv 7\ (\bmod\ 16)\ge 23$.
It desuspends seven dimensions whose Hopf invariant is
$\sigma^2_{2n-15}$.
\end{thm}
By abuse of notation, $\sigma_n$  represents a generator of
$\pi^n_{n+7}$ and $\pi_{n+7}(\S^n)$ for $n\ge 8$, respectively.
Combining the results of \cite{M1}, \cite{M2}, \cite{T}, Corollary \ref{CC}.(3), Proposition \ref{Gnp0}
and Theorem \ref{Mah},
\begin{multline}\label{wsigma}
\sharp[\iota_n,\sigma_n]=\left\{\begin{array}{ll}
240,&\ \mbox{if} \ n\ \mbox{is even and}\ n\ge 10;\\
2,&\ \mbox{if} \ \mbox{$n$ is odd and $n\ge 9$ unless $n=11$}\\
& \ \mbox{or $n\equiv 15$}\ (\bmod\ 16);\\
1,&\ \mbox{if} \ n=11 \ \mbox{or} \ n\equiv 15\ (\bmod\ 16).
\end{array}
\right.
\end{multline}
Whence, by means of Lemma \ref{lead}, the groups
$G_{n+7}(\S^n)$ for $n\ge 9$ have been fully described as well.
\section{Proof of Theorem \ref{neq}, part I}\setcounter{equation}{0}
Denote by $i_n(\mathbb{R}) : SO(n-1)\hookrightarrow SO(n)$ and
$p_n(\mathbb{R}) : SO(n)\to\S^{n-1}$
the inclusion and projection maps, respectively.
Hereafter, we use the following exact sequence induced from the fibration
$SO(n+1)\stackrel{SO(n)}{\longrightarrow}\S^n$:
$$
(\mathcal{SO}^n_k) \hspace{3mm} \pi_{k+1}(\S^n)\rarrow{\Delta}\pi_k(SO(n))
\rarrow{i_*}\pi_k(SO(n+1))\rarrow{p_*}\pi_k(\S^n)
\longrightarrow\cdots,
$$
where $i=i_{n+1}(\R)$ and $p=p_{n+1}(\R)$.
\par Since $SO(n)\cong
SO(n-1)\times\S^{n-1} \
\mbox{for} \ n=4,8,$ we get that %\addtocounter{equation}{-2}
\begin{equation} \label{37}
\Delta\pi_{k+1}(\S^n)=0, \ \mbox{if} \ n=3,7.
\end{equation}
By the exact sequence $(\mathcal{SO}^n_n)$ and the fact that
$\pi_n(SO(n))\cong\Z$ for $n\equiv 3\ (\bmod\ 4)$ \cite{K}, we have
\begin{equation} \label{Deta3}
\Delta\eta_n=0, \ \mbox{if} \ n\equiv 3\ (\bmod\ 4).
\end{equation}
We recall the formula \cite[Lemma 1]{K}
\begin{equation}\label{DeE}
\Delta(\alpha\circ E\beta) = \Delta\alpha\circ\beta.
\end{equation}
By (\ref{Deta3}) and (\ref{DeE}),
\begin{equation} \label{Deta32}
\Delta(\eta^2_n)=0, \ \mbox{if} \ n\equiv 3\ (\bmod\ 4).
\end{equation}
Denote by $V_{n,k}$ the Stiefel manifold consisting of $k$-frames in $\R^n$ for $k\le n-1$.
Then, we show
\begin{lem}\k \label{50} {\em (1)} $\Delta(\nu^2_n)=0, \ \mbox{if} \ n\equiv 5\
(\bmod\ 8);$

{\em (2)} $\Delta(\nu^2_{4n})=0$, if $n$ is odd.
\end{lem}
{\bf Proof.} Since $\pi_7(SO(5))\cong\Z$ \cite{K},
$\Delta : \pi_8(\S^5)\to\pi_7(SO(5))$ is trivial and $\Delta\nu_5=0$. So, by
(\ref{DeE}), $\Delta(\nu^2_5)=0$. Let now $n\equiv 5\ (\bmod\ 8)\ge 13$.
We consider the exact sequence $(\mathcal{SO}^n_{n+5})$:
$$
\pi_{n+6}(\S^n)\stackrel{\Delta}{\to}\pi_{n+5}(SO(n))
\stackrel{i_*}{\to}\pi_{n+5}(SO(n+1))\to 0.
$$
By \cite{B-M}, we obtain
$$\pi_{n+5}(SO(n))\cong\pi_{n+5}(SO)\oplus\pi_{n+6}(V_{n+8,8}).$$
In the light of \cite{H-M}, $\pi_{n+6}(V_{n+8,8})\cong\Z_8$ and by \cite{Bott},
$\pi_{n+5}(SO)=0$. So, $\pi_{n+5}(SO(n))\cong\Z_8$.
By \cite{K}, $\pi_{n+5}(SO(n+1))\cong\Z_8$.
From the fact that $\pi_{n+6}(\S^n)=\{\nu^2_n\}\cong\Z_2$, we
obtain $\Delta(\nu^2_n)=0$, and hence (1) follows.
\par We obtain
$\pi_9(SO(4))\cong\pi_9(SO(3))\oplus\pi_9(\S^3)\cong(\Z_3)^2$,
and so $\Delta(\nu^2_4)=0$. Let now $n\ge 3$. Then, we consider the
exact sequence
$(\mathcal{SO}^{4n}_{4n+5}$):
$$
\pi_{4n+6}(\S^{4n})\stackrel{\Delta}{\to}\pi_{4n+5}(SO(4n))
\stackrel{i_*}{\to}\pi_{4n+5}(SO(4n+1))\to 0.
$$
By \cite{K}, $\pi_{4n+5}(SO(4n+1))\cong\Z_2$.
By \cite{KM}, ${i_{13}(\R)}_*: \pi_{17}(SO(12))\to\pi_{17}(SO(13))\cong\Z_2$
is an isomorphism, and hence $\Delta(\nu^2_{12})=0$.
Let $n$ be odd and $n\ge 5$. In the light of
\cite{B-M},
$$
\pi_{4n+5}(SO(4n))\cong\pi_{4n+5}(SO)\oplus\pi_{4n+6}(V_{4n+8,8}), \
\mbox{if} \ n\ge 4.
$$
By means of \cite{Bott} and \cite{H-M},
$
\pi_{4n+5}(SO)\cong\Z_2$ and
$\pi_{4n+6}(V_{4n+8,8}) = 0$.
Hence, we obtain $\Delta(\nu^2_{4n}) = 0$ if $n$ is odd with
$n\ge 5$. This leads to (2) and completes the proof.
\hfill$\square$

\cite[Theorem 1.3]{K-M} suggests the non-triviality of
$[\iota_n,\nu^2_n]$ for $n\equiv 0,1,2,$ $3,6\ (\bmod\ 8)\ge 6$ and
\cite[Proposition 3.4]{No} gives an explicit proof of its non-triviality
for $n\equiv 2\ (\bmod\ 4)\ge 6$.
\par By Lemma \ref{pro}.(1) and (\ref{wnu}), we have $[\iota_n,\nu^2_n]=0$ if $n\equiv 7\ (\bmod\ 8)$ or $n=2^i - 3$ for $i\ge 3$.
In virtue of Lemma \ref{50} and by (\ref{P'}), we get that
%\addtocounter{equation}{-2}
\begin{equation} \label{5}[\iota_n,\nu^2_n] = 0, \ \mbox{if} \ n\equiv 5\
(\bmod\ 8)
\end{equation}

\noindent
and

\begin{equation} \label{4}
[\iota_n, \nu^2_n] = 0, \ \mbox{if} \ n\equiv 4\ (\bmod\ 8).
\end{equation}

Given elements $\alpha\in\pi_{n+k}(\S^n)$ and $\beta\in\pi_{n+k}(SO(n+1))$
satisfying $p_{n+1}(\R)\beta=\alpha$, then $\beta$ is called a lift of
$\alpha$ and we write $$\beta=[\alpha].$$
For $k\le n-1$, set $i_{k,n}=i_n(\R)\circ\cdots\circ i_{k+1}(\R)$.
\par Let now $n\equiv 0 \; (\bmod\ 4)\ge 8$. By \cite{B-M}, \cite{Bott} and
\cite{H-M}, $\pi_{2n+3}(SO(2n-2))\cong\Z\oplus\Z_4$. In the exact sequence
$(\mathcal{SO}^{2n-3}_{2n+3})$, $p_{2n-2}(\R)_\ast: \pi_{2n+3}(SO(2n-2))
\to\pi_{2n+3}(\S^{2n-3} )$ is an epimorphism by Lemma \ref{50}.(1). So, the
direct summand $\mathbb{Z}_4$ of $\pi_{2n+3}(SO(2n-2))$ is generated by
$[\nu^2_{2n-3}]$. By \cite{K}, $\pi_{2n+3}(SO(2n+1))\cong\Z\oplus\Z_2$ and
$\pi_{2n+3}(SO(2n+2))\cong\Z$. It follows from $(\mathcal{SO}^{2n+1}_{2n+3})$ that
the direct summand $\Z_2$ of $\pi_{2n+3}(SO(2n+1))$ is generated by
$\Delta\nu_{2n+1}$. By \cite{K}, $\pi_{2n+3}(SO(2n+k-1))\cong\Z\oplus\Z_2$
for $0\le k\le 2$. Hence, by use of $(\mathcal{SO}^{2n+k-1}_{2n+3})$
for $-1\le k\le 2$, $(i_{2n-2,2n+1})_*: \pi_{2n+3}(SO(2n-2))
\to \pi_{2n+3}(SO(2n+1))$ is an epimorphism and we get the relation
$$(i_{2n-2,2n+1})_\ast[\nu^2_{2n-3}]=\Delta\nu_{2n+1}.$$
Thus, we conclude
\begin{lem}\k \label{26}
$E^3J[\nu^2_{2n-3}]=[\iota_{2n+1},\nu_{2n+1}]$, if $n\equiv 0\ (\bmod\ 4)\ge 8$.
\end{lem}

Next, we need
\begin{lem}\k \label{Deps0}
Let $n\equiv 3\ (\bmod\ 4)\ge 7$. Then,

{\em (1)} $\{\Delta\iota_n,\eta_{n-1},2\iota_n\}=0$;

{\em (2)} $\Delta(E\{\eta_{n-1},2\iota_n,\alpha\})=0$, where
$\alpha\in\pi_k(\S^n)$ is an element satisfying $2\iota_n\circ\alpha=0$.
\end{lem}
{\bf Proof.} By the properties of Toda brackets and the fact that $2\pi_{n+1}(SO(n+1))=0$,
if $n\equiv 3\ (\bmod\ 4)\ge 7$ \cite{K}, we obtain
$$
i_{n+1}(\R)\circ\{\Delta\iota_n,\eta_{n-1},2\iota_n\}
=-\{i_{n+1}(\R),\Delta\iota_n,\eta_{n-1}\}\circ 2\iota_{n+1}
\subset 2\pi_{n+1}(SO(n+1))=0.
$$
It follows from $(\mathcal{SO}^n_{n+1})$ and (\ref{Deta32}) that
${i_{n+1}(\R)}_*: \pi_{n+1}(SO(n))\to\pi_{n+1}(SO(n+1))$ is a monomorphism.
This leads to (1).
\par By (\ref{DeE}) and (1), for any $\beta\in\{\eta_{n-1},2\iota_n,\alpha\}$, we obtain
$$
\Delta(E\beta)\in\Delta\iota_n\circ\{\eta_{n-1},2\iota_n,\alpha\}
=-\{\Delta\iota_n,\eta_{n-1},2\iota_n\}\circ E\alpha=0.
$$
This leads to (2) and completes the proof.
\hfill$\square$

\par We recall that $\varepsilon_{n-1}\in\{\eta_{n-1},2\iota_n,\nu^2_n\}$ and
$\mu_{n-1}\in\{\eta_{n-1},2\iota_n,E^{n-5}\sigma'''\}$
for $n\ge 5$. So, by (\ref{37}) and Lemma \ref{Deps0}.(2), we obtain

\begin{exam}\k \label{Deps}
Let $n\equiv 3\ (\bmod\ 4)$. Then,

{\em (1)} $\Delta\varepsilon_n=0$;

{\em (2)} $\Delta\mu_n=0$.
\end{exam}

\bigskip

\par Hereafter, we use often the EHP sequence of the following type:
$$(\mathcal{PE}^n_{n+k}) \hspace{3ex}
\pi^{2n+1}_{n+k+2}\rarrow{P}\pi^n_{n+k}\rarrow{E}\pi^{n+1}_{n+k+1}.$$

It is well-known that
$$
H[\iota_n,\iota_n]=0 \ \mbox{for} \ n \ \mbox{odd, and} \
H[\iota_n,\iota_n]=\pm 2\iota_{2n-1} \ \mbox{for} \ n \ \mbox{even}.
$$
So, by \cite[Proposition 2.5]{T}, we obtain
\begin{lem} \label{HP}
Suppose that $\pi^{2n+1}_{n+k+2}=E^3\pi^{2n-2}_{n+k-1}$ and $\pi^{2n-1}_{n+k}
=E\pi^{2n-2}_{n+k-1}$. Then $HP(\pi^{2n+1}_{n+k+2})=0$ for $n$ odd and
$HP(\pi^{2n+1}_{n+k+2})=2\pi^{2n-1}_{n+k}$ for $n$ even.
\end{lem}

Suppose that $\Delta\alpha = 0$ for $\alpha\in\pi_k(\S^{n-1})$. Then, by \cite{WG0}, we obtain
\begin{equation} \label{HJ}
H(J[\alpha])=\pm E^n\alpha.
\end{equation}

Now, we show
\bigskip

\noindent
    {\bf I.\  \mbox{\boldmath $[\iota_n,\nu^2_n]\ne 0$} if
    \mbox{\boldmath $n\equiv 1 (\bmod\ 8)\ge 9$}.}

    \bigskip

\par In virtue of \cite[Theorem 10.3]{T} and its proof,
$[\iota_9,\nu^2_9]=\bar{\nu}_9\nu^2_{17}\ne 0$.
\par Let $n\equiv 0\ (\bmod\ 4)\ge 8$.
Assume that $E^3(J[\nu^2_{2n-3}]\circ\nu_{4n+1})
=0\in\pi^{2n+1}_{4n+7}$. Then, by use of $(\mathcal{PE}^{2n}_{4n+6})$
and the fact that $(i_{2n-2,2n})_\ast[\nu^2_{2n-3}]$ generates the direct summand
$\Z_2$ of $\pi_{2n+3}(SO(2n))$, we obtain
$E^2(J[\nu^2_{2n-3}]\circ\nu_{4n+1})=8a[\iota_{2n},\sigma_{2n}]$ for $a\in\{0,1\}$.
By means of \cite[Proposition 11.11.i)]{T}, there exists an element
$\beta\in\pi^{2n-2}_{4n+4}$ such that $P(8\sigma_{4n+1}) = E^2\beta$ and $H\beta\in\{2\iota_{4n-5},
\eta_{4n-5},8\sigma_{4n-4}\}_2$. By the properties of Toda brackets, we see
that
\begin{eqnarray*}
\{2\iota_{4n-5},\eta_{4n-5},8\sigma_{4n-4}\}_2
\subset\{2\iota_{4n-5},\eta_{4n-5},8\sigma_{4n-4}\}\subset\\
\{2\iota_{4n-5},0,4\sigma_{4n-4}\}=2\iota_{4n-5}\circ\pi^{4n-5}_{4n+4}
+\pi^{4n-5}_{4n-3}\circ 4\sigma_{4n-3}=0.
\end{eqnarray*}
So, there exists an element $\beta'\in\pi^{2n-3}_{4n+3}$ such that $\beta
=E\beta'$. Hence, $E^2(J[\nu^2_{2n-3}]\circ\nu_{4n+1})=aE^3\beta'$.
\par In virtue of Lemma \ref{pro}.(1) and (\ref{weta}), $[\iota_{2n-1},
\eta_{2n-1}\sigma_{2n}]=0$. In the light of (\ref{P'}) and Example \ref{Deps}.(1), $[\iota_{2n-1},
\varepsilon_{2n-1}]=0$, and so $P\pi^{4n-1}_{4n+7}=0$. Therefore, by $(\mathcal{PE}^{2n-1}_{4n+5})$,
\begin{equation*}
E:\pi^{2n-1}_{4n+5}\to\pi^{2n}_{4n+6} \ \mbox{is a
monomorphism if} \ n\equiv 0\ (\bmod\ 4)\ge 8.
\end{equation*}
Hence, $E(J[\nu^2_{2n-3}]\circ\nu_{4n+1})=aE^2\beta'$.
By use of $(\mathcal{PE}^{2n-2}_{4n+4})$ and (\ref{HJ}), we have a contradictory relation
$\nu^3_{4n-5} = 0$. Thus, we get $[\iota_{2n+1},\nu^2_{2n+1}]
=E^3(J[\nu^2_{2n-3}]\circ\nu_{4n+1})\ne 0$.

\bigskip

\par We note the relation
\begin{equation} \label{pDeio}
p_n(\R)(\Delta\iota_n)=(1+(-1)^n)\iota_{n-1}.
\end{equation}
Let $n\equiv 0\ (\bmod\ 8)\ge 8$.
By use of $(\mathcal{SO}^{n-1}_{n+1})$ and \cite{K}, we get that
$i_n(\R)_*: \pi_{n+1}(SO(n-1))\to\pi_{n+1}(SO(n))$ is a monomorphism. So, we obtain
\begin{equation}\label{Denu7}
\Delta\nu_{n-1}=0,\ \mbox{if} \ n\equiv 0\ (\bmod\ 8)\ge 8.
\end{equation}
Hence, by Lemma \ref{50}.(2), $\nu_{n-1}$ and $\nu^2_{n-4}$
are lifted to $[\nu_{n-1}]\in\pi_{n+2}(SO(n))$ and $[\nu^2_{n-4}]
\in\pi_{n+2}(SO(n-3))$, respectively. We show the following

\begin{lem}\k \label{KT}
Let $n\equiv 0\ (\bmod\ 8)\ge 16$. Then, for some odd $x$,
$$
2[\nu_{n-1}]-\Delta\nu_n=x(i_{n-3,n})_*[\nu^2_{n-4}].
$$
\end{lem}
{\bf Proof.}
By use of $(\mathcal{SO}^{n-k}_{n+2})$ for $2\le k\le 4$, Lemma
\ref{50} and \cite{K}, we see that
$(i_{n-3,n-1})_*:\pi_{n+2}(SO(n-3))\to\pi_{n+2}(SO(n-1))\cong\Z_8$ is
an isomorphism and $\pi_{n+2}(SO(n-3))=\{[\nu^2_{n-4}]\}$.
In virtue of \cite{K}, $\pi_{n+2}(SO(n+1))\cong\Z_8$ and $\pi_{n+2}(SO(n))\cong\Z_{24}
\oplus\Z_8$. So, by $(\mathcal{SO}^{n-k}_{n+2})$ for $k=0,1$, we get
$\pi_{n+2}(SO(n))=\{\Delta\nu_n,[\nu_{n-1}]\}$.
By (\ref{pDeio}), we obtain $p_n(\R)(\Delta\nu_n)=2\nu_{n-1}$,
and hence $2[\nu_{n-1}]-\Delta\nu_n\in\Im\ (i_n(\R)_*:\pi_{n+2}(SO(n-1))
\to\pi_{n+2}(SO(n)))$, where $i_n(\R)_*$ is a split monomorphism. Since
$\sharp( 2[\nu_{n-1}]-\Delta\nu_n)=8$, we have the required relation and this
completes the proof.
\hfill$\square$

The relation in \cite[Lemma 11.17]{T} is regarded as the $J$-image
of that in Lemma \ref{KT}.

Now, we present a proof of the non-triviality of $[\iota_n, \nu^2_n]$ in
the case $n\equiv 0\ (\bmod\ 8)\ge 8$.

\bigskip

\noindent
{\bf II.\ $\mbox{\boldmath $[\iota_n,\nu^2_n]\ne 0$}$ if
$\mbox{\boldmath $n\equiv 0\ (\bmod\ 8)\ge 8$}$.}

\bigskip

By \cite[(7.19), Theorem 7.7]{T}, $[\iota_8,\nu^2_8]=\nu_8\sigma_{11}\nu_{18}\ne 0$. Let $n\equiv 0 \ (\bmod\ 8)\ge 16$.
In the light of \cite{B-M}, \cite{Bott} and \cite{H-M}, $\pi_{n+5}(SO(n))\cong(\Z_2)^2$.
So, by (\ref{DeE}) and Lemma \ref{KT},
$$
\Delta(\nu^2_n)=(i_{n-3,n})_*([\nu^2_{n-4}]\nu_{n+2})
$$
and hence $[\iota_n,\nu^2_n]=E^3(J[\nu^2_{n-4}]\circ\nu_{2n-1})$.
\par Assume that $E^3(J[\nu^2_{n-4}]\circ\nu_{2n-1})=0$. Then,
$E^2(J[\nu^2_{n-4}]\circ\nu_{2n-1})\in P\pi^{2n-1}_{2n+6}=\{[\iota_{n-1},\sigma_{n-1}]\}$.
By \cite[Prposition 11.11.ii)]{T}, it holds $\pi^{2n-3}_{2n+5}\subset E^2\pi^{n-4}_{2n+1}$ and
in virtue of (\ref{des7}), we have $E^2(J[\nu^2_{n-4}]\circ\nu_{2n-1})=aE^7(\gamma\sigma_{2n-10})$ for $a\in\{0,1\}$.
So, by using $(\mathcal{PE}^{n-2-k}_{2n+3-k})$ for $k=0,1$, we get that
\begin{eqnarray*}
J[\nu^2_{n-4}]\circ\nu_{2n-1}-aE^5(\gamma\sigma_{2n-10})-E\beta
\in P\pi^{2n-5}_{2n+4}
\end{eqnarray*}
for some $\beta\in\pi^{n-4}_{2n+1}$.
Hence, Lemma \ref{HP} and (\ref{HJ}) imply a contradictory relation $\nu^3_{2n-7}=0$, and thus
$[\iota_n,\nu^2_n]\ne 0$.

\vspace{2mm}

\par We note that Nomura has a different proof from II.

\section{Proof of Theroem \ref{neq}, part II}\setcounter{equation}{0}
\par Let $\omega_n(\mathbb{R})\in\pi_{n-1}(O(n))$, $\omega_n(\C)
\in\pi_{2n}(U(n))$ and $\omega_n(\H)\in\pi_{4n+2}(Sp(n))$ be the characteristic
elements for the orthogonal $O(n)$, unitary $U(n)$ and symplectic $Sp(n)$
groups, respectively. We note that $\omega_n(\R)=\Delta\iota_n$ and %\label{orDeio}
$\sharp(\Delta\iota_n)=2 \ \mbox{for odd} \ n\ge 9$.
\par Let $r_n: U(n)\to SO(2n)$ and $c_n: Sp(n)\to SU(2n)$ be the canonical maps,
respectively. Set $i_n(\C): U(n-1)\hookrightarrow U(n)$ for the inclusion map.
As it is well-known,
$$
i_{2n+1}(\R)r_n\omega_n(\C)=\omega_{2n+1}(\mathbb{R})
\;\;\mbox{and}\;\;
i_{2n+1}(\C)c_n\omega_n(\H)=\omega_{2n+1}(\C).
$$
Let $$\tau'_{2n}=r_n\omega_n(\C)\in\pi_{2n}(SO(2n))\;\mbox{and}\;
\bar{\tau}'_{4n}=r_{2n}c_n\omega_n(\H)\in\pi_{4n+2}(SO(4n)).$$
It is well-known %\cite{J}
that
$$p_{2n}(\R)\tau'_{2n}=(n-1)\eta_{2n-1}\;\mbox{and}\;
{p_{4n}(\R)}\bar{\tau}'_{4n}=\pm(n+1)\nu_{4n-1}\ \mbox{for}\ n\ge 2.$$
Whence, we obtain

\begin{lem} \label{ap}
\mbox{\em (1)} If $n$ is even and $n\ge 4$ then $i_{n+1}(\R)\tau'_n
=\Delta\iota_{n+1}$ and $p_n(\R)\tau'_n=(\frac{n}{2}-1)\eta_{n-1}$;

\mbox{\em (2)} If $n\equiv 0\ (\bmod\ 4)\ge 8$ then $(i_{n,n+2})\bar{\tau}'_n
=\tau'_{n+2}$ and $p_n(\R)\bar{\tau}'_n=\pm(\frac{n}{4}+1)\nu_{n-1}$.
\end{lem}
By use of $(\mathcal{SO}^{4n+1}_{4n+2})$, Lemma \ref{ap} and
\cite{K}, we obtain %\addtocounter{equation}{-7}
\begin{equation}\label{fr0}
\Delta(\eta^2_{4n+1})=4i_{4n+1}(\R)\bar{\tau}'_{4n}, \;\mbox{if}\, n\ge 2.
\end{equation}
So, by $(\mathcal{SO}^{4n}_{4n+2})$, we have $\tau'_{4n}\eta^2_{4n}-4\bar{\tau}'_{4n}=a\Delta\nu_{4n}$
for $a\in\{0,1,\ldots,23\}$. Composing $p_{4n}(\R)$ with this relation, using the relation
$\eta^3_{4n-1}=12\nu_{4n-1}$, (\ref{DeE}), (\ref{pDeio}) and Lemma \ref{ap},
$a$ is even and
\begin{equation}\label{fr}
\tau'_{4n}\eta^2_{4n}\equiv 4\bar{\tau}'_{4n}\ (\bmod\ 2\Delta\nu_{4n}),\;\mbox{if}\;
n\ge 2.
\end{equation}
\par Set $\tau_{2n}=J\tau'_{2n}\in\pi_{4n}(\S^{2n})\;\mbox{and}\;
\bar{\tau}_{4n}=J\bar{\tau}'_{4n}\in\pi_{8n+2}(\S^{4n})$. Then, we note that
\begin{equation} \label{tau1}
E\tau_{2n}=[\iota_{2n+1},\iota_{2n+1}], H\tau_{2n}=(n-1)\eta_{4n-1}
\end{equation}
and
\begin{equation} \label{tau2}
E^3\bar{\tau}_{4n}=[\iota_{4n+3},\iota_{4n+3}],
H\bar{\tau}_{4n}=\pm(n+1)\nu_{8n-1}
\end{equation}
By (\ref{fr0}), we have
\begin{equation}\label{4Ebt}[\iota_{4n+1},\eta^2_{4n+1}]
=4E\bar{\tau}_{4n}.\end{equation}
\par  Let $\iota_X$ be the identity class of a
space $X$. Denote by $\P^n(2)$ the Moore space of type $(\Z_2,n-1)$ and by
$i_n:\S^{n-1}\hookrightarrow\P^n(2)$, $p_n:\P^n(2)\to\S^n$ the inclusion and
collapsing maps, respectively. We recall from \cite{T1} that
\begin{equation} \label{ietap}
2\iota_{\P^n(2)}=i_n\eta_{n-1}p_n, \;\mbox{if}\; n\ge 3.
\end{equation}
Let $\bar{\eta}_n\in[\P^{n+2}(2),\S^n]\cong\Z_4$
and \mbox{$\tilde{\eta}_n\in\pi_{n+2}(\P^{n+1}(2))\cong\Z_4$} for $n\ge 3$
be an extension and a coextension of $\eta_n$, respectively.
We note that
\begin{equation} \label{exteta}
\bar{\eta}_n\in\{\eta_n,2\iota_{n+1},p_{n+1}\}, \;\mbox{if}\; n\ge 3
\end{equation}
and
\begin{equation}\label{tildeta}
\tilde{\eta}_n\in\{i_{n+1},2\iota_n,\eta_n\}, \;\mbox{if}\; n\ge 3.
\end{equation}
We have
\begin{equation} \label{ieta2p}
2\bar{\eta}_n=\eta^2_np_{n+2}\;\;\mbox{and}\;\;
2\tilde{\eta}_n=i_{n+1}\eta^2_n,\; \mbox{if}\; n\ge 3.
\end{equation}
We recall that $\bar{\eta}_n\tilde{\eta}_{n+1}=\pm E^{n-3}\nu'$ for
$n\ge 3$. Furthermore, we recall that $\pi_{n+8}(\S^n)=\{\varepsilon_n\}\cong\Z_2$ for
$3\le n\le 5$ and $\varepsilon_3=\{\eta_3,E\nu',\nu_7\}$. We need
%\addtocounter{thm}{2}
\begin{lem}\ \label{eps}
$\varepsilon_n=\{\eta_n\bar{\eta}_{n+1}, \tilde{\eta}_{n+2}, \nu_{n+4}\}$ if
$n\ge 7$.
\end{lem}
{\bf Proof.} By (\ref{tildeta}), we obtain
$$
\tilde{\eta}_7\circ\nu_9\in\{i_8, 2\iota_7, \eta_7\}\circ\nu_9
= -(i_8\circ\{2\iota_7, \eta_7, \nu_8\})\subset i_8\circ\pi_{12}(\S^7)=0.
$$
So, we can take
$$
\varepsilon_5=\{\eta_5, 2\nu_6, \nu_9\}=\{\eta_5, \bar{\eta}_6\tilde{\eta}_7,
\nu_9\}=\{\eta_5\bar{\eta}_6, \tilde{\eta}_7, \nu_9\}
$$
and
$$
\varepsilon_n\in\{\eta_n\bar{\eta}_{n+1}, \tilde{\eta}_{n+2}, \nu_{n+4}\}, \
\mbox{if} \ n\ge 5.
$$
The indeterminacy of this bracket is
$$\eta_n\bar{\eta}_{n+1}\circ\pi_{n+8}(\P^{n+3}(2))+\pi_{n+5}(\S^n)\circ\nu_
{n+5}.$$
For $n\ge 5$, by use of the homotopy exact sequence of a pair
$(\P^{n+3}(2),\S^{n+2})$, we obtain $\pi_{n+8}(\P^{n+3}(2))=\{i_{n+3}\nu^2_{n+2}\}$,
and so $\bar{\eta}_{n+1}\circ\pi_{n+8}(\P^{n+3}(2))=\{\eta_{n+1}\nu^2_{n+2}\}=0.$
Hence, the indeterminacy is trivial for $n\ge 7$ and this completes the proof.

\hfill$\square$

By \cite{K}, $\pi_{4n}(SO(4n))\cong(\Z_2)^3 \ \mbox{or} \ (\Z_2)^2,\ \mbox{if}\ n\ge 2.$ So,
\begin{equation} \label{ordtau}
\sharp\tau'_{4n}=2,\ \mbox{if}\ n\ge 2.
\end{equation}
Now, we prove
\begin{lem}\ \label{hard}
$\tau'_{n-1}\eta_{n-1}\varepsilon_n\equiv 0\,(\bmod\ \tau'_{n-1}\nu^3_{n-1})$
if $n\equiv 1\ (\bmod\ 8)\ge 17$.
\end{lem}
{\bf Proof.}
By (\ref{DeE}), (\ref{fr}), (\ref{ieta2p}) and the relation
$\nu_{n-1}\eta_{n+2}=0$ for $n\ge 7$, we have
$\tau'_{n-1}\eta^2_{n-1}\bar{\eta}_{n+1}=0$ if $n\equiv 1\ (\bmod\
8)\ge 9$.
In virtue of Lemma \ref{eps}, we obtain
\begin{eqnarray*}
\tau'_{n-1}\eta_{n-1}\varepsilon_n
=\tau'_{n-1}\eta_{n-1}\circ\{\eta_n\bar{\eta}_{n+1}, \tilde{\eta}_{n+2},
\nu_{n+4}\}\\
=-\{\tau'_{n-1}\eta_{n-1}, \eta_n\bar{\eta}_{n+1}, \tilde{\eta}_{n+2}\}
\circ\nu_{n+5}.
\end{eqnarray*}
By (\ref{ordtau}) and noting that $2\iota_{SO(n-1)}\circ\alpha=2\alpha$ for any
$\alpha\in\pi_{n+5}(SO(n-1))$, we obtain
$$
2\{\tau'_{n-1}\eta_{n-1},\eta_n\bar{\eta}_{n+1},\tilde{\eta}_{n+2}\}
=-\{2\iota_{SO(n-1)}, \tau'_{n-1}\eta_{n-1},\eta_n\bar{\eta}_{n+1}\}
\circ\tilde{\eta}_{n+3}.
$$
We have
$$\{2\iota_{SO(n-1)}, \tau'_{n-1}\eta_{n-1},\eta_n\bar{\eta}_{n+1}\}
\subset[\P^{n+4}(2),SO(n-1)].$$
By \cite{K}, $\pi_{n+3}(SO(n-1))=0$ if $n\equiv 1\ (\bmod\ 8)\ge 9$. So,
$[\P^{n+4}(2),SO(n-1)]=\pi_{n+4}(SO(n-1))\circ p_{n+4}$. In virtue of
$(\mathcal{SO}^{n-k}_{n+4})$ for $k=1,2,$ \cite{B-M}, \cite{Bott}, \cite{H-M}
%$$
%\pi_{n+5}(S^{n-1})\rarrow{\Delta}\pi_{n+4}(SO(n-1))\rarrow{i_*}\pi_{n+4}(SO(n))
%\rarrow{}0
%$$
and (\ref{Denu7}), we see that
$$
\pi_{n+4}(SO(n-1))=\{\Delta(\nu_{n-1}),[\nu_{n-2}]\}\circ\nu_{n+1}\cong(\Z_2)^2.
$$
Hence, we obtain $[\P^{n+4}(2),SO(n-1)]\circ\tilde{\eta}_{n+3}=
\pi_{n+4}(SO(n-1))\circ\eta_{n+4}=0$. This leads to the relation
$$2\{\tau'_{n-1}\eta_{n-1},\eta_n\bar{\eta}_{n+1},\tilde{\eta}_{n+2}\}=0.$$
By \cite{B-M}, \cite{Bott} and \cite{H-M}, we see that
$$
\pi_{n+5}(SO(n))\cong\Z_{16}\oplus\Z_2, \ \mbox{if}\ n\equiv 1\
(\bmod\ 8)\ge 17,$$
where the direct summand $\Z_2$ is generated by
$\Delta(\nu^2_n)$
and
\begin{equation}\label{new1}
\pi_{n+5}(SO(n-1))\cong\left\{\begin{array}{ll}
(\Z_{16})^2\oplus\Z_2\oplus\Z_{15},&\mbox{if} \ n\equiv 1\ (\bmod\ 16)\ge 17;\\
\Z_{32}\oplus\Z_8\oplus\Z_2\oplus\Z_{15},&\mbox{if} \ n\equiv 9\ (\bmod\ 16)
\ge 25,\end{array}\right.
\end{equation}
where the direct summand $\Z_2$ is generated by $\tau'_{n-1}\nu^2_{n-1}$.
%Here, it is easily seen that the first direct summand $\Z_2$ is generated by
%$\Delta(\nu^2_n)$ and the rest direct summand $\Z_2$ is generated by
%$\tau'_{n-1}\nu^2_{n-1}$.
Thus, for $n\equiv 1\ (\bmod\ 16)\ge 17$, we obtain
$$\{\tau'_{n-1}\eta_{n-1},\eta_n\bar{\eta}_{n+1},\tilde{\eta}_{n+2}\}
\subset\{\tau'_{n-1}\nu^2_{n-1}\}+8\pi_{n+5}(SO(n-1)),$$
and so $\tau'_{n-1}\eta_{n-1}\varepsilon_n\equiv 0\
(\bmod\ \tau'_{n-1}\nu^3_{n-1})$.
For $n\equiv 9\ (\bmod\ 16)\ge 25$, by $(\mathcal{SO}^{n-1}_{n+5})$,
we get
$$
\pi_{n+5}(SO(n-1);2)=\{\beta,\Delta\sigma_{n-1}-2\beta,\tau'_{n-1}\nu^2_{n-1}\}\cong
\Z_{32}\oplus\Z_8\oplus\Z_2,$$
where $\beta$ is such an element that $16\beta=8\Delta\sigma_{n-1}$ and
$i_n(\R)\beta$ is a generator of the direct summand $\Z_{16}$ of
$\pi_{n+5}(SO(n))$. Then, by (\ref{DeE}) and the relation $\sigma_{n-1}\nu_{n+6}=0$ for
$n\ge 13$ \cite[(7.20)]{T}, we have
$$
4(\Delta\sigma_{n-1}-2\beta)\circ\nu_{n+5}
=4\Delta(\sigma_{n-1}\nu_{n+6})+\beta\circ 8\nu_{n+5}=0.
$$
Consequently, we obtain the relation $\tau'_{n-1}\eta_{n-1}\varepsilon_n
\equiv 0\,(\bmod\ \tau'_{n-1}\nu^3_{n-1})$ if $n\equiv 9\ (\bmod\ 16)\ge
25$ and this completes the proof.
\hfill$\square$

\par Next, we show
\begin{lem}\k \label{n48}
If $n\equiv 0, 1\ (\bmod\ 4)\ge 8$ then $[\iota_n,\alpha]\ne 0$ for
$\alpha=\varepsilon_n, \bar{\nu}_n$ and $\eta_n\sigma_{n+1}$.
\end{lem}
{\bf Proof.}
We show $[\iota_n,\varepsilon_n]\ne 0$. Let $n\equiv 0 \ (\bmod\ 4)\ge 8$.
By \cite[Proposition 11.10.i)]{T}, there exists an element $\beta\in\pi^{n-1}_{2n+6}$ such that
$E\beta=[\iota_n,\varepsilon_n]$ and $H\beta=\eta_{2n-3}\varepsilon_{2n-2}$.
Assume that $[\iota_n,\varepsilon_n]=0$. Then, by $(\mathcal{PE}^{n-1}_{2n+6})$,
we have $\beta\in P\pi^{2n-1}_{2n+8}$.
This induces a contradictory relation
$\eta_{2n-3}\varepsilon_{2n-2}=0$, and hence
$[\iota_n,\varepsilon_n]\not=0$. Next, consider the case $n\equiv 1\ (\bmod\ 4)
\ge 9$. Then, by (\ref{tau1}),
$[\iota_n,\varepsilon_n]=E(\tau_{n-1}\varepsilon_{2n-2})$ and
$H(\tau_{n-1}\varepsilon_{2n-2})=\eta_{2n-3}\varepsilon_{2n-2}$.
Assume that $[\iota_n,\varepsilon_n]=0$.
Then, $(\mathcal{PE}^{n-1}_{2n+6})$ and Lemma \ref{HP} lead to a contradictory relation
$\eta_{2n-3}\varepsilon_{2n-2}=0$, and so $[\iota_n,\varepsilon_n]\ne 0$.
For other elements, the argument goes ahead similarly.
\hfill$\square$

\bigskip
By (\ref{JDel}) and Lemma \ref{n48}, $\Delta: \pi_{n+8}(\S^n)\to\pi_{n+7}(SO(n))$
is a monomorphism, if $n\equiv 0,1\ (\bmod\ 4)\ge 12$. So, by
$(\mathcal{SO}^n_{n+8})$, we obtain the exact sequence
\begin{multline} \label{exact1}
\pi_{n+9}(\S^n)\rarrow{\Delta}\pi_{n+8}(SO(n))\rarrow{i_*}\pi_{n+8}(SO(n+1))
\rarrow{}0,\\
\mbox{if} \ n\equiv 0,1\ (\bmod\ 4)\ge 12.\\
\end{multline}
We recall from \cite{M2} and \cite{T} the following:
\begin{equation} \label{myu}
\sharp[\iota_n,\mu_n]=\left\{\begin{array}{ll}
1,&\ \mbox{if} \ n=2,6\ \mbox{or} \ n\equiv 3\ (\bmod\ 4);\\
2,&\ \mbox{if} \ n\equiv 0,1,2\ (\bmod\ 4)\ge 4\ \mbox{unless} \ n=6,
\end{array}\right.
\end{equation}

\begin{equation} \label{emyu}
\sharp[\iota_n,\eta_n\mu_{n+1}]=\left\{\begin{array}{ll}
1,& \ \mbox{if} \ n=5 \ \mbox{or} \ n\equiv 2, 3\ (\bmod\ 4);\\
2,&\ \mbox{if} \ n\equiv 0, 1\ (\bmod\ 4)\ge 4 \ \mbox{unless} \ n=5,
\end{array}\right.
\end{equation}

\begin{equation} \label{eta2s}
\sharp[\iota_n,\eta^2_n\sigma_{n+2}]=\left\{\begin{array}{ll}
1,&\ \mbox{if} \ n\equiv 2,3\ (\bmod\ 4)\ge 6;\\
2,&\ \mbox{if} \ n\equiv 0\ (\bmod\ 4)\ge 8
\end{array}\right.
\end{equation}
and
$$(\ast) \hspace{7ex} \sharp[\iota_n,\eta^2_n\sigma_{n+2}]=2, \ \mbox{if} \
n\equiv 1\ (\bmod\ 8)\ge 17.$$

Now, we conclude
\begin{cor}\k  \label{hard2}
{\em (1)} $[\iota_n,\eta_n\varepsilon_{n+1}]\equiv 0\,(\bmod\
[\iota_n,\nu^3_n])$
if $n\equiv 1\ (\bmod\ 8)\ge 9$;

{\em (2)} $[\iota_n,\nu^3_n]=0$ if $n\equiv 5\ (\bmod\ 8)$ and
$[\iota_n,\eta_n\varepsilon_{n+1}]=[\iota_n,\eta^2_n\sigma_{n+2}]=0$
provided $n\equiv 5\ (\bmod\ 8)\ge 13$\ unless $n\equiv 53\ (\bmod\ 64)$.
\end{cor}
{\bf Proof.} We have $[\iota_9,\eta_9\varepsilon_{10}]
=\eta_9\sigma_{10}\eta_{17}\varepsilon_{18}+\sigma_9\eta^2_{16}\varepsilon_{18}
=\eta^2_9\sigma_{11}\varepsilon_{18}+4\sigma_9\nu_{16}\sigma_{19}$ $=0$.
For the case $n\equiv 1\ (\bmod\ 8)\ge 17$, (1) is a direct consequence of Lemma \ref{hard}.
\par By (\ref{DeE}) and Lemma \ref{50}.(1), we have
$\Delta(\nu^3_n)=0$. So, the first assertion of (2) holds. In the light of \cite{Mimura},
the second assertion of (2) holds for $n=13$.
Let $n\equiv 5\ (\bmod\ 8)\ge 21$. We consider the exact sequence (\ref{exact1}).
By \cite{B-M}, \cite{Bott} and \cite{H-M}, we see that
\begin{multline}\label{7SO}
\pi_{n+8}(SO(n+1))\cong\left
\{\begin{array}{ll}
\Z_4\oplus\Z_2,&\ \mbox{if} \ n\equiv 5\ (\bmod\ 32)\ge 37;\\
(\Z_4)^2,&\ \mbox{if} \ n\equiv 21\ (\bmod\ 32);\\
\Z_4,&\ \mbox{if} \ n\equiv 13\ (\bmod\ 16)
\end{array}\right.
\end{multline}
and
$$\pi_{n+8}(SO(n))\cong\left \{\begin{array}{ll}
\Z_4\oplus(\Z_2)^2,& \ \mbox{if} \ n\equiv 5\ (\bmod\ 32)\ge 37;\\
(\Z_4)^2\oplus\Z_2,& \ \mbox{if} \ n\equiv 21\ (\bmod\ 64);\\
\Z_8\oplus\Z_4\oplus\Z_2,& \ \mbox{if} \ n\equiv 53\ (\bmod\ 64);\\
\Z_4\oplus\Z_2,& \ \mbox{if} \ n\equiv 13\ (\bmod\ 16).
\end{array}\right.$$
By (\ref{myu}) and (\ref{emyu}), $[\iota_n,\mu_n]\ne 0$ and $[\iota_n,\mu_n]\eta_{2n+8}\ne 0$.
So, by the group structures of $\pi_{n+8}(SO(n+k))$ for $k=0,1$, we get that
$\Delta\mu_n$ is taken as a generator of the direct summand $\Z_2$ of
$\pi_{n+8}(SO(n))$. By (\ref{DeE}) and (\ref{fr0}), we obtain
$$\Delta(\eta^2_n\sigma_{n+2})=4i_n(\R)\bar{\tau}'_{n-1}\sigma_{n+1}$$
and hence
$$\Delta(\eta^2_n\sigma_{n+2})=\left\{\begin{array}{ll}
0,& \ \mbox{if} \ n\not\equiv 53\ (\bmod\ 64);\\
4i_n(\R)\bar{\tau}'_{n-1}\sigma_{n+1}\ne 0,& \ \mbox{if} \ n\equiv 53\
(\bmod\ 64).
\end{array}\right.$$
This leads to the second assertion of (2) and the proof is complete.
\hfill$\square$

By \cite[Proposition 4.2]{Oda},
$$
[\iota_n,\eta_n\varepsilon_{n+1}]\equiv 0\,
(\bmod\ [\iota_n,\eta^2_n\sigma_{n+2}]), \ \mbox{if} \ n\equiv 1\ (\bmod\
8)\ge 9.
$$
Thus,%\addtocounter{equation}{3}
\begin{equation}\label{eeps}
[\iota_n,\eta_n\varepsilon_{n+1}]=0, \quad \mbox{if} \quad n\equiv 1\ (\bmod\ 8)
\ge 9.\end{equation}

Next, we prove
\begin{lem}\k \label{NK}
Let $n\equiv 1\ (\bmod\ 4)\ge 5$. Then
$E(\bar{\tau}_{2n-2}\nu^2_{4n-2})=[\iota_{2n-1},\bar{\nu}_{2n-1}]$ if
and only if $[\iota_{2n+1},\nu^2_{2n+1}] = 0$.
\end{lem}
{\bf Proof.}
By (\ref{tau2}), $E^3(\bar{\tau}_{2n-2}\nu^2_{4n-2})
=[\iota_{2n+1},\nu^2_{2n+1}]$. This induces the necessary condition.

Now, suppose that $[\iota_{2n+1},\nu^2_{2n+1}] = 0$. Then, by $(\mathcal{PE}^{2n}_{4n+6})$,
%$$\pi^{4n+1}_{4n+8}\rarrow{P}\pi^{2n}_{4n+6}\rarrow{E}\pi^{2n+1}_{4n+7},$$
$E^2(\bar{\tau}_{2n-2}\nu^2_{4n-2})\in P\pi^{4n+1}_{4n+8}\cong\Z_{16}$.
We assume that $E^2(\bar{\tau}_{2n-2}\nu^2_{4n-2})=8aP(\sigma_{4n+1})$ for $a\in\{0,1\}$.
By \cite[Proposition 11.11.ii)]{T}, there exists an element $\beta\in\pi^{2n-2}_{4n+4}$ such that
$$
P(8\sigma_{4n+1}) =E^2\beta\quad \mbox{and}\quad
H\beta\in\{\eta_{4n-5}, 2\iota_{4n-4}, 8\sigma_{4n-4}\}_2.
$$
We recall that
\begin{eqnarray*}
\mu_{4n-5}
\in\{\eta_{4n-5}, 2\iota_{4n-4}, 8\sigma_{4n-4}\}_2 \
(\bmod\  \nu^3_{4n-5},\eta_{4n-5}\varepsilon_{4n-4}).&
\end{eqnarray*}
Thus, we get
$$
H\beta=\mu_{4n-5}+x\nu^3_{4n-5}+y\eta_{4n-5}\varepsilon_{4n-4}\ (x,y\in\{0,1\}).
$$
By using $(\mathcal{PE}^{2n-1}_{4n+5})$ and the assumption, we have
$$E(\bar{\tau}_{2n-2}\nu^2_{4n-2})-aE\beta
\in P\pi^{4n-1}_{4n+7}
= \{P(\bar{\nu}_{4n-1}), P(\varepsilon_{4n-1})\}.
$$
By (\ref{tau1}), $P(\bar{\nu}_{4n-1})
=E(\tau_{2n-2}\bar{\nu}_{4n-4})$ and
$P(\varepsilon_{4n-1})=E(\tau_{2n-2}\varepsilon_{4n-4})$. So, by
$(\mathcal{PE}^{2n-2}_{4n+4})$,
$$\bar{\tau}_{2n-2}\nu^2_{4n-2}-a\beta-b\tau_{2n-2}\bar{\nu}_{4n-4}
-c\tau_{2n-2}\varepsilon_{4n-4}\in P\pi^{4n-3}_{4n+6}\ (b,c\in\{0,1\}).
$$
Applying  $H: \pi^{2n-2}_{4n+4}\to\pi^{4n-5}_{4n+4}$ to this equation,
using (\ref{tau1}), Lemma \ref{HP} and the relation $\eta_{4n-5}\bar{\nu}_{4n-4}=\nu^3_{4n-5}$,
we obtain $$\nu^3_{4n-5}+a(\mu_{4n-5}+x\nu^3{4n-5}+y\eta_{4n-5}\varepsilon_{4n-4})+b\nu^3_{4n-5}+c\eta_{4n-5}\varepsilon_{4n-4}=0.$$
\par By the group structure of $\pi^{4n-5}_{4n+4}$, $a=c=0$ and $b=1$, and so
$E(\bar{\tau}_{2n-2}\nu^2_{4n-2})=E(\tau_{2n-2}\bar{\nu}_{4n-4})$.
Whence the proof is complete.
\hfill$\square$

Since $\nu_n\eta_{n+3} = 0$ and $\bar{\nu}_n\eta_{n+8}=\nu^3_n$ for $n\ge 6$,
Lemma \ref{NK} implies

\begin{cor}\k \label{NK2}
If $[\iota_{8n+3}, \nu^2_{8n+3}] = 0$, then $[\iota_{8n+1},\nu^3_{8n+1}]=0$.
\end{cor}

Now, we show
\bigskip

\noindent
{\bf III.\ $\mbox{\boldmath $[\iota_n,\nu^2_n] = 0$}$
if $\mbox{\boldmath $n = 2^i - 5\ (i\ge 4)$}$.}

\bigskip

\par We recall the Mahowald element $\eta'_i\in\pi^S_{2^i}(\S^0)$ for
$i\ge 3$ \cite{M3}. We set $\eta'_{i-1,m}=\eta'_{i-1}$ on $\S^m$ for $m=2^{i-1}-2$
with $i\ge 4$, that is, $\eta'_{i-1,m}\in\pi_{2^{i-1}+m}(\S^m)$. It satisfies
the relation $H(\eta'_{i-1,m})=\nu_{2m-1}$. Then, the assertion follows directly from \cite{Ba}.

\bigskip

Finally, we show
\bigskip

\noindent
{\bf  IV.\ $\mbox{\boldmath $[\iota_n,\nu^2_n]\ne 0$}$ if $\mbox{\boldmath $n\equiv 3\ (\bmod\ 8)\ge 19$}$
unless $\mbox{\boldmath $n=2^i-5$}$.}

\bigskip

\par By III and Corollary \ref{NK2}, we obtain
$$
[\iota_n,\nu^3_n]=0, \quad \mbox{if} \quad n=2^i-7\ (i\ge 4).
$$
Hence, from (\ref{eeps}) and the relation $\eta^2_n\sigma_{n+2}=\nu_n^3+\eta_n\varepsilon_{n+1}$,
$$[\iota_n,\eta^2_n\sigma_{n+2}]=0,\ \mbox{if}\ n=2^i-7\,(i\ge 4).$$
\par Let $n\equiv 1\ (\bmod\ 8)\ge 17$. Considering the exact sequence
(\ref{exact1}), in virtue of \cite{B-M}, \cite{Bott} and \cite{H-M}, we obtain
$$
\pi_{n+8}(SO(n))\cong\Z_2\oplus\Z_2\oplus\Z_8\quad\mbox{and}\quad
\pi_{n+8}(SO(n+1))\cong\Z_2\oplus\Z_4.
$$
By (\ref{myu}), (\ref{emyu}) and ($\ast$), we know that
$$
[\iota_n,\mu_n]\ne 0, [\iota_n,\mu_n]\eta_{2n+8}\ne 0
$$
and
$$
[\iota_n,\eta^2_n\sigma_{n+2}]\ne 0.
$$
So, by (\ref{4Ebt}), we get the relation
$$
4E(\bar{\tau}_{n-1}\sigma_{2n})=[\iota_n,\eta^2_n\sigma_{n+2}]\ne 0.
$$
By $(\ast)$ and (\ref{eeps}), we obtain
$$[\iota_n, \nu^3_n]=[\iota_n,\eta^2_n\sigma_{n+2}]\ne 0,\ \mbox{if}\ n\equiv 1\
(\bmod\ 8)\ge 17 \ \mbox{and}\ n\ne 2^i-7.$$
Thus, by Corollary \ref{NK2}, we obtain the assertion.

\bigskip

\par We are in a position to assert that Mahowald's result \cite{M2}
should be stated as follows.
\begin{thm}\k
Let $n\equiv 1\ (\bmod\ 8)\ge 9$. Then $[\iota_n,\eta^2_n\sigma_{n+2}]\ne 0$ if
and only if $n\ne 2^i-7$.
\end{thm}

\section{Proof of $[\iota_{16s+7},\sigma_{16s+7}]\ne 0$ for $s\ge 1$}\setcounter{equation}{0}
\par We give a proof of the first part of Theorem \ref{Mah}.
First of all, let $n\equiv 2\ (\bmod\ 4)\ge 10$. Then, by use of $(\mathcal{SO}^n_n)$,
%fibration $SO(n+1)\stackrel{SO(n)}{\longrightarrow}\S^n$,
Lemma \ref{ap}.(1) and \cite{K}, we obtain $\pi_n(SO(n))=\{\tau'_n\}\cong\Z_4$ and
%\addtocounter{equation}{-6}
\begin{equation} \label{deta}
2\tau'_n=\Delta\eta_n, \ \mbox{if}\ n\equiv 2\ (\bmod\ 4)\ge 10.
\end{equation}
%So, by (\ref{DeE}), $\Delta(\eta_n\sigma_{n+1})=2(\tau'_n\sigma_n)$.
\par We recall from \cite[p.95-6]{T} the construction of the element $\kappa_7\in\pi_{21}(\S^7)$.
It is a representative of a Toda bracket
$$
\{\nu_7,E\alpha,E^2\beta\}_1,
$$
where $\alpha=\bar{\eta}_9\in[\P^{11}(2),\S^9]$  is an extension of $\eta_9$
and $\beta=\widetilde{\bar{\nu}}_9\in\pi_{18}(\P^{10}(2))$ is a coextension of
$\bar{\nu}_9$ satisfying $\alpha\circ E\beta=0$. Furthermore, $\kappa_n=E^{n-7}\kappa_7$
for $n\ge 7$ and set $\widetilde{\bar{\nu}}_n=E^{n-9}\widetilde{\bar{\nu}}_9$ for $n\ge 9$.
Then, we can take
$$
\kappa_n\in\{\nu_n,\bar{\eta}_{n+3},\widetilde{\bar{\nu}}_{n+4}\} \;\mbox{for}\; n\ge 7.
$$
By \cite{K}, $\pi_{n+4}(SO(n+k))\cong\Z\oplus\Z_2$ for $k=1,2$ if
$n\equiv 7\ (\bmod\ 8)$. And, by $(\mathcal{SO}^{n+2}_{n+4})$, the direct summand
$\Z_2$ of $\pi_{n+4}(SO(n+2))$ is generated by $\Delta\nu_{n+2}$. So,
the non-triviality of $[\nu_n]\eta_{n+3}\in\pi_{n+4}(SO(n+1))$ induces the
relation
${i_{n+2}(\R)}_*([\nu_n]\eta_{n+3})=\Delta\nu_{n+2}$. Because of the fact that
$[\iota_{n+2},\nu^2_{n+2}]\ne 0$, this induces a contradictory relation
$0=\Delta\nu^2_{n+2}\ne 0$. Hence, we obtain
$$
[\nu_n]\eta_{n+3}=0, \ \mbox{if} \ n\equiv 7\ (\bmod\ 8).
$$
Next, by \cite{K},
$$
\{[\nu_n],\eta_{n+3},2\iota_{n+4}\}\subset\pi_{n+5}(SO(n+1))=0, \; \mbox{if} \;
n\equiv 7\ (\bmod\ 8).
$$
So, by (\ref{exteta}), we have
$[\nu_n]\bar{\eta}_{n+3}\in\{[\nu_n],\eta_{n+3},
2\iota_{n+4}\}\circ p_{n+5}=0$ and hence we can
define a lift of $\kappa_n$ for $n\equiv 7\ (\bmod\ 8)$, as follows:
$$
[\kappa_n]\in\{[\nu_n],\bar{\eta}_{n+3},\widetilde{\bar{\nu}}_{n+4}\}\subset\pi_{n+14}(SO(n+1))
\;\mbox{for}\; n\equiv 7\ (\bmod\ 8).$$

Let $n\equiv 7\ (\bmod\ 8)\ge 15$.
By use of $(\mathcal{SO}^{n-4}_{n-4})$, $(\mathcal{SO}^{n-k}_{n-3})$ for
$k=2,3,5$, $(\mathcal{SO}^{n-l}_{n-2})$ for $2\le l\le 5$ and \cite{K}, we obtain
$$
\pi_{n-4}(SO(n-4))=\{\beta\}\cong\Z; \
\pi_{n-3}(SO(n-4))=\{[\eta^2_{n-5}]\}\cong\Z_2;
$$
$$
\pi_{n-2}(SO(n-4))=\{[\eta^2_{n-5}]\eta_{n-3},\Delta\nu_{n-4}\}\cong(\Z_2)^2;
$$
$$
\pi_{n-4}(SO(n-3))=\{i_{n-3}(\R)\beta,\Delta\iota_{n-3}\}\cong(\Z)^2;
$$
$$
\pi_{n-3}(SO(n-3))=\{[\eta_{n-4}],\Delta\eta_{n-3}\}\cong(\Z_2
)^2;
$$
$$
\pi_{n-2}(SO(n-3))=\{[\eta_{n-4}]\eta_{n-4},\Delta\eta^2_{n-3}\}\cong(\Z_2)^2;
$$
$$\pi_{n-2}(SO(n-2))=\{\Delta\eta_{n-2}\}\cong(\Z_2)^2,$$
where $\beta$ is a generator of $\pi_{n-4}(SO(n-4))$ and%\addtocounter{equation}{-17}
\begin{equation} \label{Delta2}
\Delta\eta_{n-3}
=i_{n-3}(\R)[\eta^2_{n-5}].
\end{equation}
We need
\begin{equation} \label{id}
\{p_n(\R),i_n(\R),\Delta\iota_{n-1}\}\ni\iota_{n-1} \ (\bmod \ 2\iota_{n-1}) \
\mbox{for} \ n\ge 9.
\end{equation}
Since $\sharp[\eta_{n-4}]=2$ for $n\equiv 7\ (\bmod\ 8)$, $[\eta_{n-4}]$ is lifted to
$\overline{[\eta_{n-4}]}\in[\P^{n-2}(2),SO(n-4)]$.
Since $p_{n-3}(\R)\beta=0$, we obtain
\begin{equation} \label{beeta}
\beta\eta_{n-4}=0\in\pi_{n-3}(SO(n-4)).
\end{equation}
So, by (\ref{Delta2}) and (\ref{id}), we get
\begin{multline}\label{lifter}
[\eta_{n-4}]
\in\{i_{n-3}(\R),\Delta\iota_{n-4},\eta_{n-5}\}\,(\bmod \ i_{n-3}(\R)\circ\pi_{n-3}(SO(n-4))\\
+\pi_{n-4}(SO(n-3))\circ\eta_{n-4}{}=\{\Delta\eta_{n-3}\}) \;\mbox{for}\;n\equiv 7\ (\bmod\ 8)\ge 15.
\end{multline}
By the same reason as (\ref{37}), we obtain
$\Delta(\bar{\eta}_3)=0\in[\P^4(2),SO(3)]$.
Let $n\equiv 7\ (\bmod\ 8)\ge 15$. Then, by Lemma \ref{Deps0}.(1) and (\ref{exteta}), we obtain
$$
\Delta(\bar{\eta}_{n-4})=\Delta\iota_{n-4}\circ\bar{\eta}_{n-5}
\in-\{\Delta\iota_{n-4},\eta_{n-5},2\iota_{n-4}\}\circ p_{n-3}=0.
$$
So, $\bar{\eta}_{n-4}$ is lifted to
$[\bar{\eta}_{n-4}]\in[\P^{n-2}(2),SO(n-3)]$
for $n\equiv 7\ (\bmod\ 8)$.
We note $[\bar{\eta}_{n-4}]\in\{i_{n-3}(\R),\Delta\iota_{n-4},\bar{\eta}_{n-5}\}$ for
$n\equiv 7\ (\bmod\ 8)\ge 15$. We show

\begin{lem}\k \label{extlift}
Let $n\equiv 7\ (\bmod\ 8)\ge 15$.

\mbox{\em (1)} $\overline{[\eta_{n-4}]}\in\{i_{n-3}(\R),\Delta\iota_{n-4},
\bar{\eta}_{n-5}\} \ (\bmod\ \{\Delta\bar{\eta}_{n-3}\}+\pi_{n-2}(SO(n-3))\circ p_{n-2}+K)$,
where\ $K={i_{n-3}(\R)}_*[\P^{n-2}(2),SO(n-4)]+\pi_{n-4}(SO(n-3))
\circ\bar{\eta}_{n-4}$.

\mbox{\em (2)} ${i_{n-2}(\R)}_*K\subset \{(\Delta\eta_{n-2})p_{n-2}\}$.
\end{lem}
{\bf Proof.}
By use of the cofiber sequence
$\S^{n-3}\,\rarrow{i_{n-2}}\,\P^{n-2}(2)\,\rarrow{p_{n-2}}\,\S^{n-2}$, (\ref{ietap})
and (\ref{lifter}), we get that $\overline{[\eta_{n-4}]}\in\{i_{n-3}(\R),\Delta\iota_{n-4},
\bar{\eta}_{n-5}\} \ (\bmod\ \{\Delta\bar{\eta}_{n-3}\}+\pi_{n-2}(SO(n-3))\circ p_{n-2}+K)$
and that $[\P^{n-2}(2),SO(n-4)]=\{\overline{[\eta^2_{n-5}]},$ $(\Delta\nu_{n-4})p_{n-2}\}\cong\Z_4\oplus\Z_2$,
where $\overline{[\eta^2_{n-5}]}$ is an extension of $[\eta^2_{n-5}]$ and $2\overline{[\eta^2_{n-5}]}
=[\eta^2_{n-5}]\eta_{n-3}p_{n-2}$. Hence, by (\ref{Delta2}), we see that
$$
{i_{n-4,n-2}}_*\overline{[\eta^2_{n-5}]}\in i_{n-2}(\R)\circ\{\Delta\eta_{n-3},
2\iota_{n-3},p_{n-3}\}=$$ $$-\{i_{n-2}(\R),\Delta\eta_{n-3},2\iota_{n-3}\}
\circ p_{n-2}.
$$
Since $\{i_{n-2}(\R),\Delta\eta_{n-3},2\iota_{n-3}\}\subset\pi_{n-2}(SO(n-2))
=\{\Delta\eta_{n-2}\}$, we have ${i_{n-4,n-2}}_*[\P^{n-2}(2),SO(n-4)]\subset
\{(\Delta\eta_{n-2})p_{n-2}\}$.
\par By (\ref{exteta}) and (\ref{beeta}), we have $\beta\bar{\eta}_{n-4}\in\{\beta,\eta_{n-4},
2\iota_{n-3}\}\circ p_{n-2}\subset\pi_{n-2}(SO(n-2))\circ p_{n-2}$. Hence,
we obtain ${i_{n-2}(\R)}_*(\pi_{n-4}(SO(n-3))\circ\bar{\eta}_{n-4})
\subset\{(\Delta\eta_{n-2})p_{n-2}\}$. This completes the proof.
\hfill$\square$

We show
\begin{lem}\k \label{hard3}
$(i_{n-7,n-1})_*[\kappa_{n-8}]=\Delta\bar{\nu}_{n-1}$ if $n\equiv 7\ (\bmod\ 8)\ge 15$.
\end{lem}
{\bf Proof.}
By the group structures of $\pi_{n-5}(SO(n-7+k))$ for $0\le k\le 3$ \cite{K},
we have $(i_{n-7,n-4})_*[\nu_{n-8}]=\Delta\iota_{n-4}$, and so
$$
(i_{n-7,n-1})_*[\kappa_{n-8}]\in(i_{n-4,n-1})_*\{\Delta\iota_{n-4},
\bar{\eta}_{n-5},\widetilde{\bar{\nu}}_{n-4}\}.
$$
By Lemma \ref{extlift}, we obtain
\begin{eqnarray*}
{i_{n-3}(\R)}_*\{\Delta\iota_{n-4},\bar{\eta}_{n-5},\widetilde{\bar{\nu}}_{n
-4}\}
=-\{i_{n-3}(\R),\Delta\iota_{n-4},\bar{\eta}_{n-5}\}
\circ\widetilde{\bar{\nu}}_{n-3}\equiv\\
\overline{[\eta_{n-4}]}\circ\widetilde{\bar{\nu}}_{n-3}
\in\{[\eta_{n-4}],2\iota_{n-3},\bar{\nu}_{n-3}\}
(\bmod \ [\eta_{n-4}]\circ\pi_{n+6}(\S^{n-3})\\
+\pi_{n-2}(SO(n-3))\circ\bar{\nu}_{n-2}+K\circ\widetilde{\bar{\nu}}_{n-3}).
\end{eqnarray*}
From the relation $i_{n-2}(\R)[\eta_{n-4}]=\Delta\iota_{n-2}$, we see that
\begin{eqnarray*}
(i_{n-7,n-1})_*[\kappa_{n-8}]
&\in&-i_{n-1}(\R)\circ\{\Delta\iota_{n-2},2\iota_{n-3},\bar{\nu}_{n-3}\}\\
&=&\{i_{n-1}(\R),\Delta\iota_{n-2},2\iota_{n-3}\}\circ\bar{\nu}_{n-2}.
\end{eqnarray*}
We note ${i_{n-2}(\R)}_*(K\circ\widetilde{\bar{\nu}}_{n-3})\subset
\{\Delta\eta_{n-2}\} \circ\bar{\nu}_{n-3}=\{\Delta\nu^3_{n-2}\}=0$ by Lemma \ref{extlift}
and (\ref{5}). Since $\{i_{n-1}(\R),\Delta\iota_{n-2},2\iota_{n-3}\}\equiv \Delta\iota_{n-1} \
(\bmod\ 2\Delta\iota_{n-1})$ by (\ref{id}), we have
$$\{i_{n-1}(\R),\Delta\iota_{n-2},2\iota_{n-3}\}\circ\bar{\nu}_{n-2}=
\Delta\bar{\nu}_{n-1}.
$$
This completes the proof.
\hfill$\square$

We can take $[\bar{\nu}_n]\in\{[\nu_n],\eta_{n+3},\nu_{n+4}\}$ for $n\equiv 7\ (\bmod\ 8)$.
And we obtain $2[\bar{\nu}_n]=0$ and $2[\kappa_n]\equiv[\bar{\nu}_n]\nu^2_{n+8}\ (\bmod\
[\nu_n]\zeta_{n+3})$. By using these facts and the group structures of $\pi_{n+k}(SO(n+1))$
for $k=11,12$ and $n\equiv 7 \ (\bmod\ 8)\ge 15$, we obtain
\begin{rem}{\em
A lift $[\kappa_n]\in\pi_{n+14}(SO(n+1))$ of $\kappa_n$ is taken so that
its order is two for
$n\equiv 7\ (\bmod\ 8)\ge 15$.}
\end{rem}
Let $n\equiv 2\ (\bmod \ 4)\ge 6$. By the relation $4\zeta_n=\eta^2_n\mu_{n+2}$, Lemma \ref{pro}.(1)
and (\ref{weta2}), $4[\iota_n,\zeta_n]=0$. So, by the relation
$H[\iota_n,\zeta_n] =\pm 2\zeta_{2n-1}$, we obtain
\begin{equation} \label{zeta6}
\sharp[\iota_n,\zeta_n]=4,\  \mbox{if}\  n\equiv 2\ (\bmod\ 4)\ge 6.
\end{equation}
By \cite[Proposition 4.2]{Oda}, there exists an element
$M_t\in\pi^{8t+8}_{16t+18}$ for $t\ge 0$ such that
\begin{equation} \label{des3}
[\iota_{8t+11},\iota_{8t+11}]=E^3(M_t) \ \mbox{and} \ HM_t=\nu_{16t+15}
\ (t\ge 0).
\end{equation}
Hereafter, we fix $n=16s+7\ge 23$.
By \cite{M1}, there exists a lift $[\sigma_{n-8}]\in\pi_{n-1}(SO(n-7))$ of
$\sigma_{n-8}$. By use of the exact sequences $(\mathcal{SO}^{n-k}_{n-1})$ for
$6\le k\le 8$, by the fact that $\sharp\Delta\sigma_{n-7}=240$, $\Delta\nu^2_{n-6}\ne 0$
and by (\ref{new1}), we obtain the following:
$$
\pi_{n-1}(SO(n-7))=\{[\sigma_{n-8}],\Delta\sigma_{n-7},\tau'_{n-7}\nu^2_{n-7}\}
\cong(\Z_{16})^2\oplus\Z_2\oplus\Z_{15}.
$$
By use of $(\mathcal{SO}^{n-k}_{n-1})$ for $1\le k\le 5$ and by
\cite{B-M}, \cite{Bott}, \cite{H-M} and \cite{K}, we see that $\pi_{n-1}(SO(n-k))=\{(i_{n-7,n-k})[\sigma_{n-8}]\}
\cong\Z_{4k}$ for $k=1,2,4$.
Therefore, $(i_{n-7,n-1})[\sigma_{n-8}]=\pm\tau'_{n-1}$, $(i_{n-7,n})[\sigma_{n-8}]=\Delta\iota_n$ and $\gamma$ in (\ref{des7}) is taken
as $\gamma=J[\sigma_{n-8}]$.
By (\ref{DeE}) and (\ref{deta}), $\Delta(\eta_{n-1}\sigma_n)=2(\tau'_{n-1}\sigma_{n-1})$.
Hence, by (\ref{P'}) and (\ref{deta}), $2(E^6\gamma)=[\iota_{n-1},\eta_{n-1}]$ and
$2E^6(\gamma\sigma_{2n-8})=[\iota_{n-1},\eta_{n-1}\sigma_n]$.
Assume that $E^7(\gamma\sigma_{2n-8})=[\iota_n,\sigma_n]=0$. Then, by
$(\mathcal{PE}^{n-1}_{2n+5})$ and Lemma \ref{hard3}, $E^6(\gamma\sigma_{2n-8})=
[\iota_{n-1},\bar{\nu}_{n-1}]=E^6J[\kappa_{n-7}]$.
By $(\mathcal{PE}^{n-2}_{2n+4})$, we have
$$
E^5(\gamma\sigma_{2n-8}-J[\kappa_{n-7}])\in P\pi^{2n-3}_{2n+6}.
$$
By Corollary \ref{hard2}.(2) and its proof, $P(\nu^3_{2n-3})=0$, $P\mu_{2n-3}
\ne 0$ and
$$
P(\eta^2_{2n-3}\sigma_{2n-1})=\left\{\begin{array}{ll}
0,&\mbox{if} \ n\not\equiv 53\ (\bmod\ 64);\\
4E(\bar{\tau}_{n-3}\sigma_{2n-4}),&\mbox{if} \ n\equiv 53\ (\bmod\ 64).
\end{array}\right.
$$
So, for $b$ and $x\in\{0,1\}$, we have
$$
E^5(\gamma\sigma_{2n-8}-J[\kappa_{n-7}])
=4bE(\bar{\tau}_{n-3}\sigma_{2n-4})+xP\mu_{2n-3}.
$$
By \cite[Proposition 11.10.ii)]{T}, there exists an element
$\beta\in\pi^{n-3}_{2n+3}$ such that $P\mu_{2n-3}=E\beta$ and $H\beta
=\eta_{2n-7}\mu_{2n-6}$. Then, by $(\mathcal{PE}^{n-3}_{2n+3})$, we have
$$
E^4(\gamma\sigma_{2n-8}-J[\kappa_{n-7}])
-4b\bar{\tau}_{n-3}\sigma_{2n-4}-x\beta\in P\pi^{2n-5}_{2n+5}.
$$
This induces a relation $x\eta_{2n-7}\mu_{2n-6}=0$. Hence, $x=0$ and we can set
$$
E^4(\gamma\sigma_{2n-8}-J[\kappa_{n-7}])
-4b\bar{\tau}_{n-3}\sigma_{2n-4}=yP(\eta_{2n-5}\mu_{2n-4})\;\mbox{for}\;y\in\{0,1\}.
$$
Since $H\bar{\tau}_{n-3}=\nu_{2n-7}$ and $\nu_{2n-7}\sigma_{2n-4}=0$, we have $\bar{\tau}_{n-3}\sigma_{2n-4}
=E\xi$ for an elements $\xi\in\pi^{n-4}_{2n+2}$. By \cite[Proposition 11.10.i)]{T},
there exists an element $\beta'\in\pi^{n-4}_{2n+2}$ such that
$P(\eta_{2n-5}\mu_{2n-4})=E\beta'$ and $H\beta'=\eta^2_{2n-9}\mu_{2n-7}$.
So, we have
$$
E^3(\gamma\sigma_{2n-8}-J[\kappa_{n-7}])-4b\xi-y\beta'
\in P\pi^{2n-7}_{2n+4}.
$$
This leads to a relation $y\eta^2_{2n-9}\mu_{2n-7}=0$, and hence $y=0$.
Therefore, by (\ref{des3}), we obtain ($n=16s+7$)
$$
E^3(\gamma\sigma_{2n-8}-J[\kappa_{n-7}]
-eM_{2s-1}\zeta_{2n-12})-4b\xi=0\ (e\in\{0,1\}).
$$
We consider the EHP sequence $$\pi^{n-5}_{2n+1}\stackrel{E}{\longrightarrow}
\pi^{n-4}_{2n+2}\stackrel{H}{\longrightarrow}\pi^{2n-9}_{2n+2}\stackrel{P}{\longrightarrow}\pi^{n-5}_{2n}.$$
By (\ref{zeta6}), $H\xi=4z\zeta_{2n-9}$ for $z\in\{0,1\}$, and so there exists an element
$\xi'\in\pi^{n-5}_{2n+1}$ satisfying $E\xi'=2\xi$.
Since $\pi^{2n-11}_{2n+1}=\pi^{2n-13}_{2n}=0$, there exists an element
$\xi''\in\pi^{n-7}_{2n-1}$ satisfying $E^2\xi''=\xi'$. Hence, we have
$$
E^2(\gamma\sigma_{2n-8}-J[\kappa_{n-7}]-eM_{2s-1}\zeta_{2n-12})
-2b\xi'\in P\pi^{2n-9}_{2n+3}=0,
$$
$$
E(\gamma\sigma_{2n-8}-J[\kappa_{n-7}]-eM_{2s-1}\zeta_{2n-12})
-2bE\xi''\in P\pi^{2n-11}_{2n+2}=0
$$
and
$$
\gamma\sigma_{2n-8}-J[\kappa_{n-7}]-eM_{2s-1}\zeta_{2n-12}-2b\xi''
\in P\pi^{2n-13}_{2n+1}.
$$
We note $H(M_{2s-1}\zeta_{2n-12})=\nu_{2n-15}\zeta_{2n-12}=0$. Then,
the last relation induces a contradictory relation $\sigma^2_{2n-15}
=\kappa_{2n-15}$. Thus, we obtain the non-triviality of
$[\iota_n,\sigma_n]$ if $n\equiv 7\ (\bmod\ 16)\ge 23$.

\bigskip

By Lemma \ref{hard3}, we have $[\iota_n,\bar{\nu}_n]=E^6J[\kappa_{n-7}]$
if $n\equiv 6\ (\bmod\ 8) \ge 14$. By the parallel arguments to the above, we obtain
\begin{cor}\k \label{newbar}
$[\iota_n,\bar{\nu}_n]\neq 0$, if $n\equiv 6\ (\bmod\ 8)\ge 14$.
\end{cor}

\section{Gottlieb groups of spheres with stems for $8\le k\le 13$} %{Gottlieb groups $G_{n+k}(\S^n)$ for $8\le k\le 13$}
\setcounter{equation}{0}We know that $\pi_{n+8}(\S^n)=\{\varepsilon_n\}\cong\Z_2$ for $n=4,5$ and that
$[\iota_4,\varepsilon_4]=(E\nu')\varepsilon_7\ne 0$,
$[\iota_5,\varepsilon_5]=\nu_5\eta_8\varepsilon_9\ne 0$.
It is easy to show that $G_{16}(\S^8)=\{(E\sigma')\eta_{15},\sigma_8\eta_{15}+\bar{\nu}_8
+\varepsilon_8\}\cong(\Z_2)^2$ and $G_{17}(\S^9)=\{[\iota_9,\iota_9]\}\cong\Z_2$.
So, by Lemma \ref{n48}, we get
$$
G_{n+8}(\S^n)=0, \quad \mbox{if} \quad n\equiv 0,1\ (\bmod\ 4)\ge 4 \ \mbox{unless} \ n=8,9.
$$
Let $n\equiv 3\ (\bmod\ 4)\ge 11$. Then, by Lemma \ref{pro}.(1) and (\ref{weta}),
$[\iota_n,\eta_n\sigma_{n+1}]=0$.
In virtue of (\ref{P'}) and Example \ref{Deps}.(1), we obtain
$[\iota_n,\varepsilon_n]=0$.
Thus,
$$
G_{n+8}(\S^n)=\pi_{n+8}(\S^n), \quad \mbox{if} \quad n\equiv 3\ (\bmod\ 4).
$$
\par Now, we show the following
\begin{lem}\k \label{eps2}
\mbox{\em (1)} Let $n\equiv 2\ (\bmod\ 8)\ge 10$. Then
$\Delta\varepsilon_n=0$ and $[\iota_n,\bar{\nu}_n]=[\iota_n,\eta_n\sigma_{n+1}]\ne 0$.

\mbox{\em (2)} Let $n\equiv 6\ (\bmod\ 8)\ge 14$. Then $\Delta\varepsilon_n
=2a(i_n(\R))[\nu^2_{n-2}]\nu_{n+4}$ for $a\in\{0,1,2,3\}$. And the order of
$\Delta\pi_{n+8}(\S^n)$ is four or two according as $n\equiv 22\ (\bmod\ 32)$
or not.
\end{lem}
{\bf Proof.}
Let $n\equiv 2\ (\bmod\ 4)\ge 10$. Then, by the fact that
$\pi_{n+1}(SO(n))\cong\Z$ \cite{K}, we have $\tau'_n\eta_n=0$. So, by
(\ref{DeE}), (\ref{ieta2p}) and (\ref{deta}), we obtain
$$\Delta(\eta_n\bar{\eta}_{n+1})=2\tau'_n\circ\bar{\eta}_n
=\tau'_n\circ\eta^2_np_{n+2}=0.$$
Therefore, by Lemma \ref{eps}, we get
\begin{eqnarray*}
\Delta\varepsilon_n=\Delta\iota_n\circ\varepsilon_{n-1}
&=&\Delta\iota_n\circ\{\eta_{n-1}\bar{\eta}_n,\tilde{\eta}_{n+1},\nu_{n+3}\}\\
&=&-\{\Delta\iota_n,\eta_{n-1}\bar{\eta}_n,\tilde{\eta}_{n+1}\}\circ\nu_{n+4}.
\end{eqnarray*}
We have
$$\{\Delta\iota_n,\eta_{n-1}\bar{\eta}_n,\tilde{\eta}_{n+1}\}
\subset\pi_{n+4}(SO(n)).$$
In virtue of \cite{B-M}, \cite{Bott} and \cite{H-M},
$\pi_{n+4}(SO(n))\cong\Z_{8d}$, where $d=2$ or $1$ according as
$n\equiv 2\ (\bmod\ 8)\ge 10$ or $n\equiv 6\ (\bmod\ 8)\ge 14$.
Noting the relation $4\tilde{\eta}_{n+1}=0$, we obtain
\begin{eqnarray*}
4\{\Delta\iota_n,\eta_{n-1}\bar{\eta}_n,\tilde{\eta}_{n+1}\}
&=&-\Delta\iota_n\circ\{\eta_{n-1}\bar{\eta}_n,\tilde{\eta}_{n+1},
4\iota_{n+3}\}\\
&\subset&-\Delta\iota_n\circ\pi_{n+4}(\S^{n-1})=0.
\end{eqnarray*}
This induces $\Delta\varepsilon_n\in (2d)(\pi_{n+4}(SO(n))\circ\nu_{n+4})$.
Since $4\pi_{n+7}(SO(n))=0$ by \cite{B-M}, \cite{Bott} and \cite{H-M}, we obtain the first
assertion of (1).

Let $n\equiv 6\ (\bmod\ 8)\ge 14$. By the exact sequences
$(\mathcal{SO}^{n+k}_{n+4})$
%fibrations $SO(n+k+1)\stackrel{SO(n+k)}{\longrightarrow}\S^{n+k}$
for $k=-2,-1$ and Lemma \ref{50} we get that $i_n(\R)_\ast: \pi_{n+4}(SO(n-1))
\to\pi_{n+4}(SO(n))$ is an isomorphism and
$\pi_{n+4}(SO(n-1))=\{[\nu^2_{n-2}]\}
\cong\Z_8$. This leads to the first assertion of (2).

We recall from \cite{M2} that $\sharp[\iota_n,\eta_n\sigma_{n+1}]=2$ if
$n\equiv 2\ (\bmod\ 8)\ge 10$. So, by the first half, we obtain the second
half of (1).

By (\ref{JDel}) and (\ref{wsigma}), $\Delta : \pi_{n+7}(\S^n)\to
\pi_{n+6}(SO(n))$ is a monomorphism for even $n\ge 10$.
So, by $(\mathcal{SO}^n_{n+7})$, we have the exact sequence:
$$\pi_{n+8}(\S^n)\rarrow{\Delta}\pi_{n+7}(SO(n))\rarrow{i_*}\pi_{n+7}(SO(n+1))
\rarrow{}0.$$
By \cite{B-M}, \cite{Bott} and \cite{H-M}, we know that
$$\pi_{n+7}(SO(n+1))\cong\left
\{\begin{array}{ll}
(\Z_2)^2,&\ \mbox{if} \ n\equiv 6\ (\bmod\ 16)\ge 22;\\
\Z_2,&\ \mbox{if} \ n\equiv 14\ (\bmod\ 16)\\
\end{array}\right.$$
and by (\ref{7SO}),
$$\pi_{n+7}(SO(n))\cong\left \{\begin{array}{ll}
\Z_4\oplus\Z_2,& \ \mbox{if} \ n\equiv 6\ (\bmod\ 32)\ge 38;\\
(\Z_4)^2,& \ \mbox{if} \ n\equiv 22\ (\bmod\ 32);\\
\Z_4,& \ \mbox{if} \ n\equiv 14\ (\bmod\ 16).
\end{array}\right.$$
Hence, we obtain the second half of (2). This completes the proof.
\hfill$\square$

Now, by Lemma \ref{eps2}.(1),
$$
[\iota_n,\varepsilon_n]=0 \ \mbox{and} \ [\iota_n,\bar{\nu}_n]
=[\iota_n,\eta_n\sigma_{n+1}]\ne 0, \ \mbox{if} \ n\equiv2\,(\bmod\ 8)\ge 10.$$
Whence, we conclude that
$$
G_{n+8}(\S^n)=\{\varepsilon_n\}\cong\Z_2, \
\mbox{if} \ n\equiv 2\ (\bmod\ 8)\ge 10.
$$
Next, by (\ref{JDel}) and Lemma \ref{eps2}.(2), we obtain
$$
G_{n+8}(\S^n)\ne 0, \ \mbox{if} \ n\equiv 6\,(\bmod\ 8)\ge 14 \
\mbox{unless} \ n\equiv 22\ (\bmod\ 32).
$$
By \cite{Mimura}, \cite{Oda2} and \cite{T}, we obtain $G_{14}(\S^6;2)=\pi^6_{14}$ and
$G_{n+8}(\S^n)=\{\eta_n\sigma_{n+1}\}\cong\Z_2$ if $n=14,22$.
Since $[\iota_6,[\iota_6,\alpha_1(6)]]=0$ by Proposition \ref{Toda}, we obtain
$G_{14}(\mathbb{S}^6;3)=\pi_{14}(\mathbb{S}^6;3)$. Thus, we have shown
\begin{prop}\k \label{8}
The group $G_{n+8}(\S^n)$ is equal to the following group:
$0$ if $n\equiv 0,1\ (\bmod\ 4)\ge 4$ unless $n=8,9$;
$\pi_{n+8}(\S^n)$ if $n=6$ or $n\equiv 3\ (\bmod\ 4)$;
$\{\varepsilon_n\}\cong\Z_2, \ \mbox{if} \ n\equiv 2\ (\bmod\ 8)\ge 10$.
Moreover, $G_{n+8}(\S^n)\ne 0$ if $n\equiv 6\ (\bmod\ 8)\ge 14$ unless
$n\equiv 22\ (\bmod\ 32)$, $G_{16}(\S^8)=\{(E\sigma')\eta_{15},\sigma_8\eta_{15}
+\bar{\nu}_8+\varepsilon_8\}\cong(\Z_2)^2$, $G_{17}(\S^9)=\{[\iota_9,\iota_9]\}\cong\Z_2$
and $G_{n+8}(\S^n)=\{\eta_n\sigma_{n+1}\}\cong\Z_2$ if $n=14,22$.
\end{prop}
Finally, we propose
\begin{con}\k
$[\iota_n,\eta_n\sigma_{n+1}]=0$ and $G_{n+8}(\S^n)=\{\eta_n\sigma_{n+1}\}\cong\Z_2$,
if $n\equiv 6\ (\bmod\ 8)\ge 14$.
\end{con}

\bigskip

\par Obviously, we obtain $G_{15}(\S^6)=\pi_{15}(\S^6)$ and
$G_{19}(\S^{10})=\{3[\iota_{10},\iota_{10}],$ $
\nu^3_{10},\eta_{10}\varepsilon_{11}\}\cong 3\Z\oplus(\Z_2)^2$.
Let $n\equiv 2\ (\bmod\ 4)\ge 14$. Then, by (\ref{deta}),
$$
[\iota_n,\eta^2_n\sigma_{n+2}]=[\iota_n,\eta_n\varepsilon_{n+1}]=0.
$$
By (\ref{myu}), $[\iota_n,\mu_n]\ne 0$. Whence, we obtain
$$G_{n+9}(\S^n)=\{\nu^3_n,\eta_n\varepsilon_{n+1}\}\cong(\Z_2)^2,
\; \mbox{if} \; n\equiv 2\ (\bmod\ 4)\ge 14.$$
\par Let now $n\equiv 3\ (\bmod\ 4)\ge 11$.
Then, by Lemma \ref{pro}.(1) and (\ref{weta}), $[\iota_n,\eta^2_n\sigma_{n+2}]
=[\iota_n,\eta_n\varepsilon_{n+1}]=0$ and by Example \ref{Deps}.(2), $[\iota_n,\mu_n]=0$.
Whence, we obtain
$$
G_{n+9}(\S^n)=\pi_{n+9}(\S^n), \; \mbox{if} \; n\equiv 3\ (\bmod\ 4).
$$
\par It is easily seen that $G_{13}(\S^4)=\{\nu^3_4\}\cong\Z_2$.
Let $n\equiv 4\ (\bmod\ 8)\ge 12$. By Lemma \ref{pro}.(1) and (\ref{4}), we have
$[\iota_n,\nu^3_n]=0$. In the light of (\ref{myu}) and (\ref{eta2s}),
$[\iota_n,\eta_n\varepsilon_{n+1}]=[\iota_n,\eta^2_n\sigma_{n+2}]\ne 0$
and $[\iota_n,\mu_n]\ne 0$. Assume that $P(\alpha_{2n+1}+\mu_{2n+1})=0$ for
$\alpha_{2n+1}=\eta_{2n+1}\varepsilon_{2n+2}$ or $\eta^2_{2n+1}\sigma_{2n+3}$.
By \cite[Proposition 11.10.i)]{T}, there exists an element $\beta\in\pi^{n-1}_{2n+7}$
satisfying $E\beta=0$ and $H\beta=\eta_{2n-3}(\alpha_{2n-2}+\mu_{2n-2})=\eta_{2n-3}\mu_{2n-2}$.
On the other hand, $(\mathcal{PE}^{n-1}_{2n+7})$ and Lemma \ref{HP} imply a contradictory
relation $H\beta=0$. So, $[\iota_n,\alpha_n]\ne [\iota_n,\mu_n]$ and hence
$$
G_{n+9}(\S^n)=\{\nu^3_n\}\cong\Z_2, \quad \mbox{if} \quad n\equiv 4\ (\bmod\ 8).
$$
\par Obviously, we obtain
$G_{18}(\S^9)=\{\sigma_9\eta^2_{16},\nu^3_9,\eta_9\varepsilon_{10}\}\cong(\Z_2)^3$.
Let now $n\equiv 1\ (\bmod\ 8)\ge 17$.
By (\ref{myu}), $[\iota_n,\mu_n]\ne 0$ and by (\ref{eeps}),
$[\iota_n,\eta_n\varepsilon_{n+1}]=0$. In the light of IV,
$[\iota_n,\eta^2_n\sigma_{n+2}]=0$ if $n=2^i-7$ for $i\ge 4$ and
$[\iota_n, \nu^3_n]=[\iota_n,\eta^2_n\sigma_{n+2}]\ne 0$ if
$n\equiv 1\ (\bmod\ 8)\ge 17$ and $n\ne 2^i-7$.
By the parallel argument to the case $n\equiv 4\ (\bmod\ 8)$,
we get $[\iota_n,\eta^2_n\sigma_{n+2}]\ne[\iota_n,\mu_n]$
\cite[Proposition 11.10.ii)]{T}.
So, we obtain
$$G_{n+9}(\S^n)=\left\{\begin{array}{ll}
 \{\eta_n\varepsilon_{n+1}\}\cong\Z_2, \
\mbox{if} \ n\equiv 1\ (\bmod\ 8)\ge 17 \ \mbox{and} \ n\ne 2^i-7;\\
\{\eta_n\varepsilon_{n+1},\eta^2_n\sigma_{n+2}\}\cong(\Z_2)^2, \
\mbox{if} \ n=2^i-7 (i\ge 5).
\end{array}
\right.$$
\par Obviously, we obtain $G_{14}(\S^5)=\{\nu^3_5,\eta_5\varepsilon_6\}\cong(\Z_2)^2$.
Let $n\equiv 5\ (\bmod\ 8)\ge 13$. By Corollary \ref{hard2}.(2) and
(\ref{myu}), $\nu^3_n\in G_{n+9}(\S^n)$ and $\mu_n\not\in G_{n+9}(\S^n)$.
Furthermore, by Corollary \ref{hard2}.(2), $\eta_n\varepsilon_{n+1}\in G_{n+9}(\S^n)$ unless
$n\equiv 53\ (\bmod\ 64)$. So, we obtain
\begin{eqnarray*}
G_{n+9}(\S^n)=\{\nu^3_n,\eta_n\varepsilon_{n+1}\}
\cong(\Z_2)^2, \ \mbox{if}\ n\equiv 5\ (\bmod\ 8)
\ \mbox{and} \ n\not\equiv \ 53 \ (\bmod\ 64).
\end{eqnarray*}
\par At the end, we use the following:
$$
\zeta_n\in\{2\iota_n,\eta_n,\alpha_{n+1}\}_2 \ (\bmod \ 2\zeta_n)\ \mbox{for} \
\alpha_{n+1}=\eta^2_{n+1}\sigma_{n+3}\ \mbox{or} \ \eta_{n+1}\varepsilon_{n+2}, \
\mbox{if} \ n\ge 11.
$$
Let $n\equiv 0\ (\bmod\ 8)\ge 16$. By \cite[Proposition 11.11.i)]{T}, there exists an element $\beta\in\pi^{n-2}_{2n+6}$
such that $[\iota_n,\alpha_n]=E^2\beta$ and $H\beta\in\
\{2\iota_{2n-5},\eta_{2n-5},$ $\alpha_{2n-4}\}_2
\ni\zeta_{2n-5}\ (\bmod\ 2\zeta_{2n-5})$. Assume that $[\iota_n,$ $\alpha_n]=0$.
Then, $(\mathcal{PE}^{n-1}_{2n+7})$ and (\ref{des7}) induce a relation
$E(\beta-aE^6(\gamma\eta_{2n-10}\mu_{2n-9}))=0$ for $a\in\{0,1\}$.
Hence, by $(\mathcal{PE}^{n-2}_{2n+6})$ and Lemma \ref{HP}, we have a contradictory
relation $\zeta_{2n-5}\in 2\pi^{2n-5}_{2n+6}$. Whence, we get that
$[\iota_n,\alpha_n]\ne 0$. In the light of (\ref{myu}) and (\ref{emyu}),
we know $[\iota_n,\mu_n]\ne 0$ and $[\iota_n,\mu_n]\eta_{2n+8}\ne 0$. This implies that
$[\iota_n,\alpha_n]\ne[\iota_n,\mu_n]$ and $[\iota_n,\nu^3_n]\ne[\iota_n,\mu_n]$.
It is easy to show that $[\iota_8,\nu^3_8]=\eta_8\bar{\varepsilon}_9$ and
$G_{17}(\S^8)=\{(E\sigma')\eta^2_{15},\sigma_8\eta^2_{15}+\nu^3_8+\eta_8\varepsilon_9\}
\cong(\Z_2)^2$.
By \cite{M-M-O} and \cite{T}, we obtain $G_{25}(\S^{16})=0$. By
\cite[4.14]{No2},
there exists an element $\tau_1\in\pi^{n-6}_{2n+2}$ such that
$$
[\iota_n,\nu^3_n]=E^6\tau_1,\ H\tau_1=\eta_{2n-13}\kappa_{2n-12}, \
\mbox{if} \ n\equiv 0\ (\bmod\ 8)\ge 16.
$$
Assume that $[\iota_n,\nu^3_n]=0$. Then, by $(\mathcal{PE}^{n-1}_{2n+7})$ and
(\ref{des7}), we have $E^5(\tau_1-aE^2(\gamma\eta_{2n-10}\mu_{2n-9}))=0$ for
$a\in\{0,1\}$. So, by $(\mathcal{PE}^{n-2}_{2n+6})$, we have
$E^4(\tau_1-aE^2(\gamma\eta_{2n-10}\mu_{2n-9}))\in P\pi^{2n-3}_{2n+8}
=\{[\iota_{n-2},\zeta_{n-2}]\}$.  By applying $H : \pi^{n-2}_{2n+6}\to\pi^{2n-5}_{2n+6}$
to this relation  and by (\ref{zeta6}),
$E^4(\tau_1-aE^2(\gamma\eta_{2n-10}\mu_{2n-9}))=0$.
By the fact that $\pi^{2n-5}_{2n+7}=\pi^{2n-7}_{2n+6}=0$, we obtain
$E^2(\tau_1-aE^2(\gamma\eta_{2n-10}\mu_{2n-9}))=0$.
Hence, by $(\mathcal{PE}^{n-5}_{2n+3})$ and (\ref{des3}), we have
$$
E(\tau_1-aE^2(\gamma\eta_{2n-10}\mu_{2n-9}))\in P\pi^{2n-9}_{2n+5}
=E^3M_t\circ\{\sigma^2_{2n-11},\kappa_{2n-11}\}$$
for $n=8t+16$. By $(\mathcal{PE}^{n-6}_{2n+2})$, we obtain
$$
\tau_1-E^2(a\gamma\eta_{2n-10}\mu_{2n-9}+bM_t\sigma^2_{2n-14}
+cM_t\kappa_{2n-14})\in P\pi^{2n-11}_{2n+4}$$
with $b,c\in\{0,1\}$.
This induces a contradictory relation $\eta_{2n-13}\kappa_{2n-12}
\in 2\pi^{2n-13}_{2n+2}$. Thus, we conclude that
$$
[\iota_n,\nu^3_n]\ne 0,\ \mbox{if} \ n\equiv 0\ (\bmod\ 8)\ge 16.
$$
Summing the above, we get
\begin{prop}\k \label{9}
The group $G_{n+9}(\S^n)$ is equal to the following group:
$\pi_{n+9}(\S^n)$ if $n=6$ or $n\equiv 3\ (\bmod\ 4)$;
$\{\nu^3_n,\eta_n\varepsilon_{n+1}\}\cong(\Z_2)^2$ if $n\equiv 2 \ (\bmod\ 4)\ge 14$,
$n=2^i-7$ for $i\ge 5$ or $n\equiv 5 \ (\bmod\ 8)$ unless $n\equiv 53\ (\bmod\ 64)$;
$\{\nu^3_n\}\cong\Z_2$ if $n\equiv 4\ (\bmod\ 8)$;
$\{\eta_n\varepsilon_{n+1}\}\cong\Z_2$ if $n\equiv 1\ (\bmod\ 8)
\ge 17$ and $n\ne 2^i-7$; $0$ if $n\equiv 0\ (\bmod\ 8)\ge 16$.
Moreover, $G_{17}(\S^8)=\{(E\sigma')\eta^2_{15},
\sigma_8\eta^2_{15}+\nu^3_8+\eta_8\varepsilon_9\}\cong(\Z_2)^2$,
$G_{18}(\S^9)=\{\sigma_9\eta^2_{16},\nu^3_9,\eta_9\varepsilon_{10}\}
\cong(\Z_2)^3$ and $G_{19}(\S^{10})=\{3[\iota_{10},\iota_{10}],
\nu^3_{10},\eta_{10}\varepsilon_{11}\}\cong 3\Z\oplus(\Z_2)^2$.
\end{prop}

By Lemma \ref{lead}, Corollary\ref{CC}.(3), Proposition \ref{Gnp0} and
(\ref{emyu}), we have determined $G_{n+10}(\S^n)$ for $n\ge 12$.
\par It is easily seen that
$$
G_{n+10}(\S^n)= \left\{ \begin{array}{ll}
\{\nu_4\sigma'+E\varepsilon',2E\varepsilon',
\alpha_1(4)\alpha_2(7),&\\
\nu_4\alpha_2(7),\nu_4\alpha'_1(7)\}, & \mbox{if $n=4$};\\
\pi_{15}(\S^5), & \mbox{if $n=5$};\\
\pi^6_{16}\oplus\Z_3\{3\beta_1(6)\}, & \mbox{if $n=6$};\\
\{\sigma_8\nu_{15},\nu_8\sigma_{11},\sigma_8\alpha_1(15)\}, & \mbox{if $n=8$};\\
\{\sigma_9\nu_{16},\beta_1(9)\}, & \mbox{if $n=9$};\\
\pi^{10}_{20}=\{\sigma_{10}\nu_{17},\eta_{10}\mu_{11}\}, & \mbox{if $n=10$};\\
\pi_{21}(\S^{11}), & \mbox{if $n=11$}.\\
\end{array}
\right.
$$
Thus, by summing up the above results, we get
\begin{prop}\k \label{ToMa2}
The group $G_{n+10}(\S^n)$ is isomorphic to one of the
following groups:
$\Z_{120}\oplus\Z_6$, $\Z_{72}\oplus\Z_2$, $\Z_{24}\oplus\Z_2$,
$\Z_{24}\oplus\Z_8$, $\Z_{24}$, $\Z_4\oplus\Z_2$, $\Z_6\oplus\Z_2$ according as
$n=4,5,6,8,9,10,11$. Furthermore, $G_{n+10}(\S^n)$ is isomorphic to the group:
$0$ if $n\equiv 0 \ (\bmod\ 4)\ge 12$; $\Z_2$ if $n\equiv 2\ (\bmod\ 4)\ge 14$;
$\Z_3$ if $n\equiv 1\ (\bmod\ 4)\ge 13$ and $\Z_6$ if $n\equiv 3\ (\bmod\ 4)\ge 15$.
\end{prop}

\par We recall that $\pi_{n+11}(\S^n;3)=\{\alpha_3(n)\}\cong\Z_3$ for $n=3,4$ and
that $\pi_{n+11}(\S^n;3)=\{\alpha'_3(n)\}\cong\Z_9$ for $n\ge 5$, where
$3\alpha'_3(n)=\alpha_3(n)$ for $n\ge 5$.
By \cite{Mimura}, \cite{M-M-O}, \cite{M-T}, \cite{T} and (\ref{zeta6}),
$\sharp[\iota_n,\zeta_n]=1,4,8,1,4,1,8,1$
according as $n=5,6,8,9,10,11,12,13$.
We easily obtain that
$$
G_{n+11}(\S^n)= \left\{ \begin{array}{ll}
\{\nu_4\sigma'\eta_{14},\nu_4\bar{\nu}_7,\nu_4\varepsilon_7,\\
2E\mu',\varepsilon_4\nu_{12},(E\nu')\varepsilon_7\}, & \mbox{if $n=4$};\\
\pi_{16}(\S^5), & \mbox{if $n=5$};\\
\{4\zeta_6,\bar{\nu}_6\nu_{14}\}, & \mbox{if $n=6$};\\
\{\bar{\nu}_8\nu_{16}\}, & \mbox{if $n=8$};\\
\pi_{20}(\S^9), & \mbox{if $n=9$};\\
4\pi_{21}^{10}, & \mbox{if $n=10$};\\
\pi_{22}(\S^{11}), & \mbox{if $n=11$};\\
\{3[\iota_{12},\iota_{12}]\}, & \mbox{if $n=12$}.\\
\end{array}
\right.
$$
\par By abuse of notations, $\zeta_n$ for $n\ge 5$ represents a generator of
the direct summands $\Z_8$ of $\pi^n_{n+11}$ and $\Z_{504}$ of
$\pi_{n+11}(\S^n)$, respectively.
By \cite{M2}, \cite{Mimura}, \cite{M-M-O}, \cite{T}, Corollary \ref{CC}.(3), Proposition
\ref{Gnp0} and (\ref{zeta6}), we obtain
$$
\sharp[\iota_n,\zeta_n] =\left\{\begin{array}{ll}
1,&\ \mbox{if} \ n\equiv 1,5,7\ (\bmod\ 8)\ge 5;\\
252,& \ \mbox{if} \ n\equiv 2\ (\bmod\ 4)\ge 6;\\
504,&\ \mbox{if} \ n\equiv 0\ (\bmod\ 4)\ge 8.
\end{array}
\right.
$$
Assume that $n\equiv 3\,(\bmod\,8)\ge 19$. Then, by \cite{B-M}, \cite{Bott} and \cite{H-M},
we obtain
$$
\pi_{n+10}(SO(n))\cong\left\{\begin{array}{ll}
(\Z_2)^3,&\ \mbox{if} \ n\equiv 3\ (\bmod\ 32)\ge 35;\\
(\Z_2)^2\oplus\Z_4,& \ \mbox{if} \ n\equiv 19\ (\bmod\ 64);\\
(\Z_2)^2\oplus\Z_8,& \ \mbox{if} \ n\equiv 51\ (\bmod\ 128);\\
(\Z_2)^2\oplus\Z_{16},& \ \mbox{if} \ n\equiv 115\ (\bmod\ 128);\\
(\Z_2)^2,&\ \mbox{if} \ n\equiv 11\ (\bmod\ 16)\ge 27
\end{array}
\right.
$$
and
$$
\pi_{n+10}(SO(n+1))\cong\left\{\begin{array}{ll}
(\Z_2)^4\oplus\Z_3,&\ \mbox{if} \ n\equiv 3\ (\bmod\ 32)\ge 35;\\
(\Z_2)^3\oplus\Z_{12},& \ \mbox{if} \ n\equiv 19\ (\bmod\ 64);\\
(\Z_2)^3\oplus\Z_{24},& \ \mbox{if} \ n\equiv 51\ (\bmod\ 64);\\
(\Z_2)^3\oplus\Z_3,&\ \mbox{if} \ n\equiv 11\ (\bmod\ 16)\ge 27.
\end{array}
\right.
$$
In the exact sequence $(\mathcal{SO}^n_{n+10})$, we get
$p_{n+1}(\R)_\ast([\eta_n]\mu_{n+1})=\eta_n\mu_{n+1}$ and $p_{n+1}(\R)_\ast[\beta_1(n)]
=\beta_1(n)$. So, $p_{n+1}(\R)_\ast$ is a split epimorphism. Whence, by the group structures
of $\pi_{n+10}(SO(n+k))$ for $k=0,1$, we obtain
$$\Delta\zeta_n=0, \; \mbox{if} \; n\not\equiv 115\ (\bmod\ 128) \;
\mbox{and} \; \Delta\zeta_n\ne 0, \; \mbox{if} \; n\equiv
115\ (\bmod\ 128).
$$
So, we have
$
[\iota_n,\zeta_n]=0, \; \mbox{if}\; n\equiv 3\,(\bmod\ 8)\ge 19 \;\mbox{and}\;
n\not\equiv 115\ (\bmod\ 128)$.
Consequently, by use of Lemma \ref{lead} and \cite{T}, the groups $G_{n+11}(\S^n)$ have
been determined if $n\ge 13$ except $n\equiv 115\ (\bmod\ 128)$.
Thus, by summing up the above results, we get
\begin{prop}\k \label{11}
The group $G_{n+11}(\S^n)$ is isomorphic to one of
the following groups: $(\Z_2)^6$, $\Z_{504}\oplus(\Z_2)^2$,
$\Z_2\oplus\Z_4$, $\Z_2$, $\Z_{504}\oplus\Z_2$,
$\Z_2$, $\Z_{504}$, $3\Z$ according as $n=4,5,6,8,9,10,11,12$.
%$G_{22}(\S^{11})=\pi_{22}(\S^{11})$
Furthermore, $G_{n+11}(\S^n)$ is isomorphic to the group: $\Z_{504}$ if $n\equiv 1,5,7\
(\bmod\ 8)\ge 13$; $\Z_2$ if $n\equiv 2\ (\bmod\ 4)\ge 14$; $0$
if $n\equiv 0\ (\bmod\ 4)\ge 16$ and $\Z_{504}$ if $n\equiv 3\ (\bmod\ 8)\ge 19$
provided $n\not\equiv 115\ (\bmod\ 128)$.
\end{prop}
We recall that $\zeta_n\in\{2\iota_n,\eta^3_n,\sigma_{n+3}\}$ for $n\ge 11$.
So, by the fact that $2\Delta\iota_n=0$ for $n$ odd, we obtain
$$\Delta\zeta_n=-\{\Delta\iota_n,2\iota_{n-1},\eta^3_{n-1}\}\circ\sigma_{n+3
}\ \mbox{for}\ n \ \mbox{odd}\ \mbox{and} \ n\ge 11.$$
We see that
$2\{\Delta\iota_n,2\iota_{n-1},\eta^3_{n-1}\}=-\Delta\iota_n
\circ\{2\iota_{n-1},\eta^3_{n-1},2\iota_{n+2}\}=0$.
Hence, by the fact that $\pi_{n+3}(SO(n))\cong\Z_{16}$ for
$n\equiv 3\ (\bmod\ 8)\ge 11$, we obtain
\begin{rem}
{\em $\Delta\zeta_n\in 8(\pi_{n+3}(SO(n))\circ\sigma_{n+3})$ if
$n\equiv 3\ (\bmod\ 8)\ge 11$.}
\end{rem}
Finally, we recall $\pi_{22}(\S^{10})=\{[\iota_{10},\nu_{10}]\}\cong\Z_{12}$.
By Proposition \ref{Toda}, $G_{22}(\S^{10})=\pi^{10}_{22}$. It is easily seen,
in the light of \cite{T}, Corollary \ref{CC}.(3) and Proposition \ref{Gnp0},
that $G_{n+12}(\S^n)=\pi_{n+12}(\S^n)$ unless $n=10$ and
$$
G_{n+13}(\S^n)= \left\{ \begin{array}{ll}
\{\nu^2_4\sigma_{10},\nu_4\eta_7\mu_8,(E\nu')\eta_7\mu_8,\\
(\alpha_1(4)+\nu_4)\beta_1(7)\}\cong\Z_{24}\oplus(\Z_2)^2, & \mbox{if $n=4$};\\
\{3[\iota_{14},\iota_{14}]\}\cong 3\Z, & \mbox{if $n=14$};\\
\pi^n_{n+13}, & \mbox{if $n$ is even}\\
& \mbox{unless $n=2,4,14$};\\
\pi_{n+13}(\S^n), & \mbox{if $n$ is odd}.
\end{array}
\right.
$$
%\newpage
We close the paper with the table of $G_{n+k}(\S^n)$ for $1\le k\le 13$ and $2\le n\le 26$:

\vspace{1mm}

\hspace{-1.7cm}
\begin{tabular}{||c||c|c|c|c|c|c|c||} \hline\hline
$G_{n+k}(\S^n)$&n=2&n=3&n=4&n=5&n=6&n=7&n=8\\ \hline\hline
k=1&$\infty$&$2$ &$0$&$0$&$2$&$2$&$0$   \\ \hline
k=2&$2$&$2$ &$0$&$2$&$2$&$2$&$0$  \\ \hline
k=3&$2$&$12$ &$3\infty+2$&$24$&$2$&$24$&$0$  \\ \hline
k=4&$12$&$2$ &$(2)^2$&$2$&$0$&$0$&$0$  \\ \hline
k=5&$2$&$2$ &$(2)^2$&$2$&$3\infty$&$0$&$0$  \\ \hline
k=6&$2$&$3$ &$24+3$&$2$&$0$&$2$&$0$  \\ \hline
k=7&$3$&$15$ &$0$&$30$&$0$&$120$&$3\infty+2$\\ \hline
k=8&$15$&$2$ &$0$&$0$&$24+2$&$(2)^3$&$(2)^2$  \\ \hline
k=9&$2$&$(2)^2$&$2$&$(2)^2$&$(2)^3$&$(2)^4$&$(2)^2$  \\ \hline
k=10&$(2)^2$&$12+2$ &$120+6$&$72+2$&$24+2$&$24+2$&$24+8$  \\ \hline
k=11&$12+2$&$84+(2)^2$ &$(2)^6$&$504+(2)^2$&$4+2$&$504+2$&$2$  \\ \hline
k=12&$84+(2)^2$&$(2)^2$ &$(2)^6$&$(2)^3$&$240$&$0$&$0$  \\ \hline
k=13&$(2)^2$&$6$ &$24+(2)^2$&$6+2$&$2$&$6$&$(2)^2$ \\ \hline\hline
\end{tabular}

%\newpage
\vspace{2mm}

\hspace{-1.7cm}
\begin{tabular}{||c||c|c|c|c|c|c|c|c|c||} \hline\hline
$G_{n+k}(\S^n)$&n=9&n=10&n=11&n=12&n=13&n=14&n=15&n=16&n=17\\ \hline\hline
k=1&$0$&$0$ &$2$&$0$&$0$&$0$&$2$&$0$&$0$   \\ \hline
k=2&$0$&$2$ &$2$&$0$&$0$&$2$&$2$&$0$&$0$  \\ \hline
k=3&$12$&$2$ &$12$&$2$&$24$&$2$&$24$&$0$&$12$  \\ \hline
k=4&$0$&$0$ &$0$&$0$&$0$&$0$&$0$&$0$&$0$  \\ \hline
k=5&$0$&$0$ &$0$&$0$&$0$&$0$&$0$&$0$&$0$  \\ \hline
k=6&$0$&$0$ &$2$&$2$&$2$&$0$&$2$&$0$&$0$  \\ \hline
k=7&$120$&$0$ &$240$&$0$&$120$&$0$&$240$&$0$&$120$\\ \hline
k=8&$2$&$2$ &$(2)^2$&$0$&$0$&$2$&$(2)^2$&$0$&$0$  \\ \hline
k=9&$(2)^3$&$3\infty+(2)^2$&$(2)^3$&$2$&$(2)^2$&$(2)^2$&$(2)^3$&$0$&$2$  \\ \hline
k=10&$24$&$4+2$ &$6+2$&$0$&$3$&$2$&$6$&$0$&$3$  \\ \hline
k=11&$504+2$&$2$ &$504$&$3\infty$&$504$&$2$&$504$&$0$&$504$  \\ \hline
k=12&$0$&$4$ &$2$&$(2)^2$&$2$&$0$&$0$&$0$&$0$  \\ \hline
k=13&$6$&$2$ &$6+2$&$(2)^2$&$6$&$3\infty$&$3$&$0$&$3$ \\ \hline\hline
\end{tabular}

\vspace{8mm}

\hspace{-1.2cm}
\begin{tabular}{||c||c|c|c|c|c|c|c|c|c||} \hline\hline
$G_{n+k}(\S^n)$&n=18&n=19&n=20&n=21&n=22&n=23&n=24&n=25&n=26\\ \hline\hline
k=1&$0$&$2$&$0$&$0$&$0$&$2$&$0$&$0$&$0$   \\ \hline
k=2&$2$&$2$&$0$&$0$&$2$&$2$&$0$&$0$&$2$  \\ \hline
k=3&$2$&$12$&$0$&$12$&$2$&$24$&$0$&$12$&$2$  \\ \hline
k=4&$0$&$0$&$0$&$0$&$0$&$0$&$0$&$0$&$0$  \\ \hline
k=5&$0$&$0$&$0$&$0$&$0$&$0$&$0$&$0$&$0$  \\ \hline
k=6&$0$&$0$&$2$&$2$&$0$&$2$&$0$&$0$&$0$  \\ \hline
k=7&$0$&$120$&$0$&$120$&$0$&$120$&$0$&$120$&$0$ \\ \hline
k=8&$2$&$(2)^2$&$0$&$0$&$2$&$(2)^2$&$0$&$0$&$2$  \\ \hline
k=9&$(2)^2$&$(2)^3$&$2$&$(2)^2$&$(2)^2$&$(2)^3$&$0$&$2$&$(2)^2$  \\ \hline
k=10&$2$&$6$&$0$&$3$&$2$&$6$&$0$&$3$&$2$  \\ \hline
k=11&$2$&$504$&$0$&$504$&$2$&$504$&$0$&$504$&$2$  \\ \hline
k=12&$0$&$0$&$0$&$0$&$0$&$0$&$0$&$0$&$0$ \\ \hline
k=13&$0$&$3$&$0$&$3$&$0$&$3$&$0$&$3$&$0$ \\ \hline\hline
\end{tabular}

\bigskip

Like in \cite{T}, an integer $n$ indicates the cyclic group $\Z_n$ of order $n$, the symbol $\infty$ an infinite
cyclic group $\Z$, the symbol $+$ the direct sum of groups and $(2)^k$ indicates the direct sum of $k$-copies
of $\Z_2$.

\vspace{5mm}
{\sc
\noindent
Faculty of Mathematics\\
and Computer Science\\
Nicholas Copernicus University\\
87-100 Toru\'n, Chopina 12/18\\
Poland}\\
e-mail: marek@mat.uni.torun.pl

\vspace{2mm}

{\sc
\noindent
Department of Mathematical Sciences\\
Faculty of Science, Shinshu University\\
Matsumoto 390-8621, Japan}\\
e-mail: jmukai0@gipac.shinshu-u.ac.jp

\end{document}